\documentclass[pdflatex,sn-mathphys-num]{sn-jnl}%

\usepackage{graphicx}%
\usepackage{amsmath}%
\usepackage{amssymb}%
\usepackage{amsfonts}%
\usepackage{amsthm}%
\usepackage{booktabs}%
\usepackage{algorithm}%
\usepackage{mathtools}%
\usepackage{enumerate}%

\newcommand{\N}{\mathbb{N}}
\newcommand{\R}{\mathbb{R}}
\newcommand{\C}{\mathbb{C}}
\newcommand{\V}{\mathbb{V}}
\newcommand{\W}{\mathbb{W}}
\newcommand{\Pro}{\mathbb{P}}

\DeclareMathOperator{\GL}{GL}
\DeclareMathOperator{\U}{U}
\DeclareMathOperator{\rk}{rk}
\newcommand{\EE}{\hspace{0.3mm}{\text{\raisebox{\depth}{\rotatebox{180}{$E$}\hspace{-1em}}}}\hspace{0.25cm}}
\newcommand\scalemath[2]{\scalebox{#1}{\mbox{\ensuremath{\displaystyle #2}}}}

\newcommand{\bz}{\mathbf{z}}
\newcommand{\bw}{\mathbf{w}}
\newcommand{\by}{\mathbf{y}}
\newcommand{\bx}{\mathbf{x}}
\newcommand{\bu}{\mathbf{u}}
\newcommand{\bv}{\mathbf{v}}
\newcommand{\ba}{\mathbf{a}}
\newcommand{\bb}{\mathbf{b}}
\newcommand{\bt}{\mathbf{t}}
\newcommand{\bs}{\mathbf{s}}
\newcommand{\bo}{\mathbf{0}}
\DeclareMathOperator{\Row}{Row}
\DeclareMathOperator{\spn}{span}
\DeclareMathOperator{\Iso}{Iso}
\DeclareMathOperator{\Id}{Id}

\DeclareMathOperator{\MV}{MV}
\DeclareMathOperator{\conv}{conv}
\DeclareMathOperator{\sing}{sing}
\DeclareMathOperator{\reg}{reg}

\DeclareMathOperator{\vHD}{vHDdeg}
\DeclareMathOperator{\HD}{HDdeg}
\DeclareMathOperator{\vH}{v\mathcal{H}}
\DeclareMathOperator{\vPH}{v\mathcal{PH}}
\DeclareMathOperator{\vX}{v\Xi}
\DeclareMathOperator{\vHDp}{vHDpol}
\DeclareMathOperator{\HDp}{HDpol}
\DeclareMathOperator{\ED}{EDdegree}
\DeclareMathOperator{\rED}{\mathbb{R}EDdegree}

\theoremstyle{thmstyleone}%
\newtheorem{Theo}{Theorem}[section]%
\newtheorem{Prop}[Theo]{Proposition}%
\newtheorem{Coro}[Theo]{Corollary}%
\newtheorem{Lemma}[Theo]{Lemma}%

\theoremstyle{thmstyletwo}%
\newtheorem{ex}{Example}[section]%
\newtheorem{rmk}{Remark}[section]%
\newtheorem{Conj}{Conjecture}[section]%

\theoremstyle{thmstylethree}%
\newtheorem{Def}{Definition}[section]%

\begin{document}

\title[Article Title]{The Hermitian Distance degree of an Algebraic Variety}

\author[1]{\fnm{Davide} \sur{Furch\`{i}}}\email{dfurchi@uninsubria.it}

\affil[1]{\orgdiv{Dipartimento di Scienza e Alta Tecnologia}, \orgname{Universit\'{a} degli Studi dell'Insubria}, \orgaddress{\city{Como}, \postcode{22100}, \country{Italy}}}

\abstract{In this paper we develop an algebraic theory to study the problem of finding the minimum distance point from an algebraic variety with respect to the Hermitian distance function. The theory generalizes the Euclidean Distance degree introduced in \cite{dhost2014}, replacing a positive symmetric bilinear form by a Hermitian form. Various examples are presented to show the robustness of the machinery.}




\maketitle

\section{Introduction}\label{sec:intro}
In this paper we consider the problem of finding the critical points of an induced distance function from an algebraic variety. In \cite{dhost2014} the case of the Euclidean distance function is studied, we adapt the framework to study a distance function induced by a Hermitian form on a complex vector space.

We use mostly geometric algebraic tools to introduce notions that characterize the problem. For any algebraic variety we define a natural number, that we call \emph{virtual Hermitian Distance degree} and we denote $\vHD$. This number describes the level of complexity in finding the nearest point from the variety and bounds from above the number of critical points. Even simple examples present nontrivial challenges.

The main difference from the Euclidean case is that this bound is rarely attained. Moreover, the number of critical points depends on the position and thus we call the set of generic values \emph{Hermitian Distance degree} and we denote it $\HD$.

We start by doing some introductory observations on the problem and by recalling the Euclidean distance degree in Section~\ref{sec:pre}. In Section~\ref{sec:first} we present the definitions for both the general and the projective case. We prove several general results, give a different characterization of the problem, compute several values of $\vHD$ and $\HD$ and discuss the notion of duality. In Section~\ref{sec:hddisc} we discuss an important object that relates to the number of critical points and to how this number changes. Moreover, we define an interesting construction that possesses similar properties to the evolute of a real curve. In Section~\ref{sec:hdpoly} we introduce a polynomial that encodes all the important information about the Hermitian distance problem. We also show some of the properties of this polynomial and discuss a few examples.

\section{Preliminaries}\label{sec:pre}

Let $\V$ be a complex vector space with complex conjugation, that is an antilinear map $\bar{\phantom{v}}\colon\V\to\V$ such that $\bar{\phantom{v}}^2=\Id$. We recall that an antilinear map $f$ is additive and conjugate homogeneous in the sense that $f(\lambda\bz)=\bar{\lambda}f(\bz)$ for any $\bz\in\V$ and $\lambda\in\C$.

Consider a differentiable function
\begin{align*}
    q\colon&\V\times\V\to\C\\
    &(\bz,\bw)\mapsto q(\bz,\bw)
\end{align*}
such that $q(\bz,\bw)=\overline{q(\bw,\bz)}$. We are interested in the point that minimizes the induced real valued function
\begin{align*}
    q_{\bu}\colon X&\subseteq\V\to\R\\
    &\bz\quad\longmapsto q(\bu-\bz,\bu-\bz)
\end{align*}
where $\bu\in\V$ and $X=V(f_1,\ldots,f_s)=V(I_X)\subseteq\V$ is an algebraic variety with $f_1,\ldots,f_s\in\C[\bz]$.\newline

We say that $X$ is a real algebraic variety if the polynomials $f_1,\ldots,f_s$ have real coefficients.

Let $\V_{\R}\subseteq\V$ be the subspace invariant under the action of the conjugation. The similar problem of considering a real algebraic variety $X$ and a symmetric bilinear form $\hat{q}$ which restriction on $\V_{\R}$ is positive was first addressed in \cite{dhost2014}. We recall here the basic notions introduced by the authors.

For $\hat{q}$ a (positive) symmetric bilinear form, if $\bz\in X_{\reg}$ is a critical point of $\hat{q}_{\bu}$, then $\nabla_{\bz}\hat{q}_{\bu}(\bz)$ is Euclidean orthogonal to the tangent space $T_{\bz}X$ of $X$ at $\bz$. In the case $\hat{q}$ is the canonical Euclidean inner product, the \emph{critical ideal} of $X$ is
\begin{equation}\label{conded}
    \left( I_{X}+\left\langle\text{$(c+1)$-minors of $\begin{bmatrix}
              \bu-\bz\\
              J(f)
  \end{bmatrix}$}\right\rangle\right)\colon\left( I_{X_{\sing}}\right)^{\infty}\subseteq\C[\bz,\bu],
\end{equation}
where $c$ is the codimension of $X$. The \emph{Euclidean distance degree} of $X$ is the constant number of solutions of the critical ideal \eqref{conded} for $\bu$ in a Zariski-open subset of $\V$, and it is denoted $\ED(X)$. However, to make a geometric sense of critical points for the Euclidean distance we need to restrict our attention to real solutions. The critical points are only the solutions in $\V_{\R}$ of the critical ideal of a point $\bu\in\V_{\R}$. To make this distinction clearer we introduce the following notation.

\begin{Def}
Let $X\subseteq\V$ be a real algebraic variety and $\hat{q}$ be a (positive) symmetric bilinear form on $\V$. The \emph{real Euclidean distance degree} of $X$ is the subset of $\N$ consisting in the numbers of critical points of $\hat{q}_{\bu}$ which do not vary for $\bu$ in an Euclidean-open subset of $\V_{\R}$, and it is denoted $\rED(X)$.
\end{Def}

By definition, it holds $\max\rED(X)\leq\ED(X)$ and it is common to have equality $\ED=\max\rED$. As an example, for the conics there hold the equalities
\begin{align*}
	&\text{$\ED(X)=3$ and $\rED(X)=\lbrace 1,3\rbrace$ if $X$ is a parabola,}\\
	&\text{$\ED(X)=2$ and $\rED(X)=\lbrace 2\rbrace$ if $X$ is a circle,}\\
	&\text{$\ED(X)=4$ and $\rED(X)=\lbrace 2,4\rbrace$ if $X$ is an ellipse or a hyperbola.}
\end{align*}\newline

We now consider the Hermitian distance problem.

The conjugation map extends on complex vector bundles over $X$. In particular, this holds true for the tangent bundle $TX$ of $X$.

In this work we will use the notations $\partial_z$ and $\partial_{\bar{z}}$ for the complex linear operators given by the Wirtinger derivatives which  satisfy the usual rules of differentiation. We will also use the notation $\nabla_{\bz}$ to denote the gradient with respect to the Wirtinger derivatives $\partial_z$ of a differentiable function.

We will use the notation $\|\phantom{z}\|_{p}$ to denote the $p$-norm for $p=1,2$. Moreover, we will denote with $\langle\phantom{z},\phantom{z}\rangle_{\R}$ the canonical Euclidean inner product and with $\perp_{\R}$ the perpendicularity condition $\langle\phantom{z},\phantom{z}\rangle_{\R}=0$.

We prove a basic lemma that characterizes the induced distance function. By some abuse of notation, we will call (complex) points that satisfy a certain perpendicularity condition \emph{critical points}. We will see in the next result and with more emphasis in Proposition~\ref{rel-com} that these points and critical points of a function with real domain are in some sense interchangeable.

\begin{Lemma}\label{perpcond}
Let $X$ be an algebraic variety and $\bu$ be a point. A regular critical point $\bz\in X$ of the function $q_{\bu}$ satisfies $\nabla_{\bz}q_{\bu}(\bz)\perp_{\R}T_{\bz}X$.
\end{Lemma}
\begin{proof}
The variety $X$ is locally a subset of $\C^n$ defined by a finite number $c\leq n$ of polynomials $f_1,\ldots f_c\in\C[\bz]$, where $c$ is the local codimension of $X$. Without loss of generality, we can reorder the $z_k$ in such a way that for a regular point 
\begin{equation*}
\bz=(\bz_a,\bz_b)=(z_1,\ldots,z_{n-c},z_{n-c+1},\ldots,z_n)\in X,
\end{equation*}
it holds $\det J_{\bz_b}(f_1,\ldots,f_c)\neq 0$ and, by the implicit function theorem, there exists a holomorphic map
\begin{equation*}
  h\colon U_1\subseteq\C^{n-c}\to U_2\subseteq\C^c
\end{equation*}
where $U_1$, $U_2$ are open sets and $U_1\times U_2$ is a neighborhood of $\bz$ such that the map 
\begin{equation*}
  (\bz_a,h(\bz_a))\colon U_1\subseteq\C^{n-c}\to U_1\times U_2\subseteq\C^n
\end{equation*}
is a parametrization of $X$ around $\bz$. Consider now the composition 
\begin{equation*}
  q_{\bu}(\bz_a,h(\bz_a))\colon U_1\subseteq\C^{n-c}\to\R.
\end{equation*}
Since it is real-valued, if we see this map on an open set of $\R^{2(n-c)}$ with coordinates 
\begin{equation*}
z_k=x_k+iy_k\qquad\text{for $k=1,\ldots,n-c$},
\end{equation*}
then its partial derivatives $\partial_{x_k}$ and $\partial_{y_k}$ for $k=1,\ldots,n-c$ vanish in a critical point. Since again this composition is real-valued, for any $k=1,\ldots,n-c$ there hold the equations
\begin{equation*}
   \partial_{x_k}q_{\bu}(\bx_a+i\by_a,h(\bx_a+i\by_a))=\partial_{y_k}q_{\bu}(\bx_a+i\by_a,h(\bx_a+i\by_a))=0
\end{equation*}
if and only if it holds the equation
\begin{equation*}
  \quad\partial_{z_k}q_{\bu}(\bz_a,h(\bz_a))=\frac{\partial_{x_k}-i\partial_{y_k}}{2}q_{\bu}(\bx_a+i\by_a,h(\bx_a+i\by_a))=0.
\end{equation*}
Thus, using the complex chain rule in a critical point we obtain
\begin{align}
  \bo=\nabla_{\bz_a}q_{\bu}(\bz_a,h(\bz_a))&=\begin{bmatrix}
      \nabla_{\bz_{a}}q_{\bu} & \nabla_{\bz_b}q_{\bu}
  \end{bmatrix}\begin{bmatrix}
      I_{n-c}\\
      J_{\bz_a}(h)
  \end{bmatrix}\notag\\
  &=\nabla_{\bz_a}q_{\bu}+\nabla_{\bz_b}q_{\bu}\cdot J_{\bz_a}(h)\label{jacobian}
\end{align}
and
\begin{align*}
  \bo&=\nabla_{\bz_a}\left(f_1(\bz_a,h(\bz_a)),\ldots,f_c(\bz_a,h(\bz_a))\right)\\
  &=\begin{bmatrix}
      J_{\bz_a}(f_1,\ldots,f_c) & J_{\bz_b}(f_1,\ldots,f_c)
  \end{bmatrix}\begin{bmatrix}
      I_{n-c}\\
      J_{\bz_a}(h)
  \end{bmatrix}\\
  &=J_{\bz_a}(f_1,\ldots,f_c)+J_{\bz_b}(f_1,\ldots,f_c)\cdot J_{\bz_a}(h)
\end{align*}
which implies $J_{\bz_a}(h)=-J_{\bz_b}(f_1,\ldots,f_c)^{-1}\cdot J_{\bz_a}(f_1,\ldots,f_c)$. Substituting this last identity in \eqref{jacobian} we obtain
\begin{equation*}
  \nabla_{\bz_a}q_{\bu}=\nabla_{\bz_b}q_{\bu}\cdot J_{\bz_b}(f_1,\ldots,f_c)^{-1}\cdot J_{\bz_a}(f_1,\ldots,f_c).
\end{equation*}
Thus, we can write
\begin{align*}
  \nabla_{\bz}q_{\bu}=\begin{bmatrix}
      \nabla_{\bz_a}q_{\bu} & \nabla_{\bz_b}q_{\bu}
  \end{bmatrix}&=\nabla_{\bz_b}q_{\bu}\cdot J_{\bz_b}(f_1,\ldots,f_c)^{-1}\cdot\begin{bmatrix}
      J_{\bz_a}(f_1,\ldots,f_c) & J_{\bz_b}(f_1,\ldots,f_c)
  \end{bmatrix}\\
  &=\nabla_{\bz_b}q_{\bu}\cdot J_{\bz_b}(f_1,\ldots,f_c)^{-1}\cdot J(f_1,\ldots,f_c).
\end{align*}
In the end, since $\nabla_{\bz}q_{\bu}\in\Row(J(f_1,\ldots,f_c))$ and the Euclidean orthogonal complement of this space is $T_{\bz}X$, a regular critical point satisfies $\nabla_{\bz}q_{\bu}(\bz)\perp_{\R}T_{\bz}X$.
\end{proof}

From now on, we consider $q$ to be a Hermitian form. Thus, we can write the orthogonality condition as $\langle\nabla_{\bz}q_{\bu}(\bz),\nu\rangle_{\R}=q(\bu-\bz,\nu)=0$ for any $\nu\in T_{\bz}X$.\newline

The Hermitian case is profoundly different from the Euclidean case. Lemma~\ref{perpcond} suggests that the conjugate operation will be very present. In particular, we will need to use the following objects.

\begin{Def}
A \emph{generalized polynomial} is an expression of the form 
\begin{equation*}
    p(\bz,\bar{\bz})=\sum_{(\alpha,\beta)\in\N^{2r}}a_{\alpha,\beta}\bz^{\alpha}\bar{\bz}^{\beta}
\end{equation*}
where the coefficients $a_{\alpha,\beta}\in\C$ are non zero only for finitely many multi-indices $(\alpha,\beta)\in\N^{2r}$. We denote its conjugate $\bar{p}(\bz,\bar{\bz})=\sum_{(\alpha,\beta)}\bar{a}_{\alpha,\beta}\bz^{\alpha}\bar{\bz}^{\beta}$.
\end{Def}

In particular we will use the \emph{Hermitian Killing form} from \cite{f2024} to count the number of solutions of systems of generalized polynomial equations.


\section{The Hermitian Distance degree}\label{sec:first}

For the sake of simplicity, we can set a basis on the complex vector space $\V$ and let $q$ be the canonical Hermitian inner product on $\V$. In this case, it holds $\nabla_{\bz}q_{\bu}(\bz)=\bar{\bz}-\bar{\bu}$ and by Lemma~\ref{perpcond} a regular critical point $\bz\in X$ satisfies $\bar{\bu}-\bar{\bz}\perp_{\R}T_{\bz}X$. We will mostly consider this situation from now on.\newline

Let $X$ be an algebraic variety, we denote $X_{\circ}\coloneqq X_{\reg}\times\overline{X}_{\reg}\subseteq\V^2$.\newline


Fix a radical ideal $I_X\coloneqq\langle f_1,\ldots,f_s\rangle\in\C[\bz]$, for our purpose we assume that $X=V(I_X)\subseteq\V$ is irreducible and
that $I_X$ is a prime ideal.

Write $J(f)$ for the $s\times n$ Jacobian matrix of $(f_1,\ldots,f_s)$. The singular locus of $X$ is defined by the ideal
\begin{equation*}
  I_{X_{\sing}}\coloneqq I_X+\langle\text{$c$-minors of $J(f)$}\rangle
\end{equation*}
where $c$ is the codimension of $X$. We now augment the matrix $J(f)$ with the row vector $\bar{\bu}-\bar{\bz}$ to get the $(s+1)\times n$-matrix
\begin{equation*}
    \begin{bmatrix}
              \bar{\bu}-\bar{\bz}\\
              J(f)
  \end{bmatrix}.
\end{equation*}
This matrix has rank $\leq c$ on the critical points of $q_{\bu}$ on $X$ since it has to hold $\nabla_{\bz}q_{\bu}(\bz)=\bar{\bz}-\bar{\bu}\in\Row(J(f))$.

We now introduce two new collections $\bw$ and $\bv$ of variables $\lbrace w_1,\ldots,w_n\rbrace$ and $\lbrace v_1,\ldots,v_n\rbrace$ respectively, where $w_k$ and $v_k$ represent the variables $\bar{z}_k$ and $\bar{u}_k$ respectively for $k=1,\ldots,n$. We extend the conjugate map to a map $\ast$ between polynomial rings
\begin{align*}
    \ast\colon\C[\bz,\bw&,\bu,\bv]\to\C[\bz,\bw,\bu,\bv]\\
    &g\qquad\mapsto\quad\ast(g)=g^{\ast}
\end{align*}
such that $z_k^{\ast}=w_k$, $w_k^{\ast}=z_k$, $u_k^{\ast}=v_k$, $v_k^{\ast}=u_k$ for $k=1,\ldots,n$ and $a^{\ast}=\bar{a}$ for any $a\in\C$. A similar map is used in \cite[Subsection 3.2]{f2024}.

Thus, using the ideal
\begin{equation*}
  I_{X}^{\prime}\coloneqq\left\langle\text{$(c+1)$-minors of $\begin{bmatrix}
              \bv-\bw\\
              J_{\bz}(f)
  \end{bmatrix}$}\right\rangle,
\end{equation*}
we define the \emph{Hermitian critical ideal} of $X$ as the following saturation
\begin{equation}\label{cond}
  \left(I_{X}+(I_{X})^{\ast}+I_{X}^{\prime}+(I_{X}^{\prime})^{\ast}\right)\colon\left( I_{X_{\sing}}\cdot (I_{X_{\sing}})^{\ast}\right)^{\infty}\subseteq\C[\bz,\bw,\bu,\bv].
\end{equation}
For fixed $\bu$ and $\bv$ we will call the ideal above, in the polynomial ring $\C[\bz,\bw]$, the Hermitian critical ideal of $(\bu,\bv)$.

Note that the Hermitian critical ideal is invariant under the action of the map $\ast$. In particular, the Hermitian critical ideal of a point $(\bu,\bar{\bu})$ is invariant under the action of the map $\ast$ restricted to the polynomial ring $\C[\bz,\bw]$.

For a more general version of the following result see Theorem~\ref{corr}.

\begin{Lemma}\label{lemmastart}
Let $X\subseteq\V$ be an algebraic variety, the variety of the Hermitian critical ideal of a generic point $(\bu,\bv)\in\V^2$ is finite. If $\bv=\bar{\bu}$ then the $\bz$ component of solutions such that $\bw=\bar{\bz}$ are exactly the critical points of the function $q_{\bu}$ on $X_{\reg}$.
\end{Lemma}
\begin{proof}
For fixed $(\bz,\bw)\in X_{\circ}$ the Jacobian  has rank $c$, so the $(c+1)\times (c+1)$-minors of $\begin{bmatrix}
    \bv-\bw\\
    J_{\bz}(f)
\end{bmatrix}$ define an affine-linear subspace of dimension $c$ in $\bv$. Similarly for $\bu$. Thus, the variety of tuples $(\bz,\bw,\bu,\bv)\in X_{\circ}\times\V^{2}$ that are solutions of the Hermitian critical ideal is irreducible of dimension $2n$. The fiber of the projection over a generic point $(\bu,\bv)\in\V^2$ must hence be finite.

The second assertion follows from the definition of critical points.
\end{proof}

The number of points in the zero locus of the Hermitian critical ideal of a point $(\bu,\bv)$ is constant on a Zariski-open subset of $\V^2$. This number counts points $(\bz,\bw)\in X_{\circ}$. The number of solutions such that $\bw=\bar{\bz}$ can be lower and if $\bv=\bar{\bu}$ it is equal to the number of critical points of $q_{\bu}$. This number is only locally constant on the projection over the $\bu$ component of the Zariski-open dense subset above and so depends on the point $\bu\in\V$. By a slight abuse of notation, we sometimes also call these \emph{critical points}, they coincide with the \emph{conjugated singles} of \cite[Definition 3.6]{f2024}. Since the definition is additive over the components of a variety, the assumption of irreducibility of $X$ is not restrictive.

\begin{Def}
Let $X$ be an algebraic variety. The \emph{virtual Hermitian Distance degree} of $X$ is the constant number of solutions of the Hermitian critical ideal \eqref{cond} of $(\bu,\bv)$ in a Zariski-open subset of $\V^2$, and it is denoted $\vHD(X)$. The \emph{Hermitian Distance degree} of $X$ is the subset of $\N$ consisting in the numbers of critical points of $q_{\bu}$ which do not vary for $\bu$ in an Euclidean-open subset of $\V$, and it is denoted $\HD(X)$.
\end{Def}

By definition, it holds the inequality $\max\HD(X)\leq\vHD(X)$. It will be common to get the strict inequality for an algebraic variety $X$, for example see the case of conics in Subsection~\ref{ssec:conics}.\newline

Some inequalities involving both the $\ED$ and $\HD$ of a real algebraic variety are presented in the following result.

\begin{Prop}\label{HDvsED}
Let $X$ be a real algebraic variety. There hold the inequalities 
\begin{equation*}
    \vHD(X)\geq\ED(X)\geq\min\HD(X)
\end{equation*}
and
\begin{equation*}
	\max\HD(X)\geq\max\rED(X).
\end{equation*}
In particular, critical points of the Euclidean distance are real critical points of the Hermitian distance.
\end{Prop}
\begin{proof}
Setting $\bv=\bu$ the left hand inequality of the chain follows from the fact that if $\bz\in X$ is in the zero locus of the critical ideal \eqref{conded} of $\bu$, since it holds $\overline{X}=X$ then $(\bz,\bz)\in X\times\overline{X}$ is in the zero locus of the Hermitian critical ideal \eqref{cond} of $(\bu,\bu)$.

For the right hand inequality of the chain note that the minimum on the right is bounded from above by the minimum on the subset in which $\bv=\bu$ is real. Now, if $(\bz,\bar{\bz})\in X\times\overline{X}$ is in the zero locus of the Hermitian critical ideal \eqref{cond} of $(\bu,\bu)$, i.e.\ is a critical point of the Hermitian distance, since all the polynomials are real and in particular $X=\overline{X}$, then $\bz,\bar{\bz}\in X$ are in the zero locus of the critical ideal \eqref{conded} of $\bu$.

The remaining inequality and the last part follow setting $\bv=\bu$ real. Now, if $\bz\in X$ real is in the zero locus of the critical ideal \eqref{conded}, i.e.\ is a real critical point of the Euclidean distance, then $(\bz,\bar{\bz})=(\bz,
\bz)\in X\times\overline{X}$ is in the zero locus of the Hermitian critical ideal \eqref{cond} of $(\bu,\bu)$, i.e.\ it is a real critical point of the Hermitian distance.
\end{proof}

All the inequalities of this last proposition could be strict in general. Moreover, in general non real solutions $\bz$ of the critical ideal \eqref{conded} of $\bu$ are not components of the solutions $(\bz,\bw)$ of the Hermitian critical ideal \eqref{cond} of $(\bu,\bu)$.

We will see in Theorem~\ref{corr} that for any non trivial algebraic variety $X$ it holds $\min\HD(X)>0$.\newline

We now briefly discuss the case of parametrized varieties.

Let $X$ be parametrized by $\psi(\bz)\colon U\subseteq\C^m\to X\subseteq\C^n$ where $U$ is an open subset. Thus, we can write
\begin{equation*}
  q_{\bu}(\bz)=\sum\limits_{k=1}^n|\psi_k(\bz)-u_k|^2
\end{equation*}
and define the map
\begin{equation*}
  D_{\bu,\bv}(\bz,\bw)\coloneqq\sum\limits_{k=1}^n(\psi_k(\bz)-u_k)(\bar{\psi}_k(\bw)-v_k).
\end{equation*}
Arguing as in Lemma~\ref{perpcond}, we can infer that regular critical points must satisfy the equality $\nabla_{\bz,\bw}D_{\bu,\bv}(\bz,\bw)=\bo$. The critical locus of $D_{\bu,\bv}(\bz,\bw)$ in the subset $U\times\overline{U}\subseteq\C^{2m}$ is the relative set of points in which the Jacobians of $\psi(\bz)$ and $\bar{\psi}(\bw)$ have maximal rank. The closure of the image of this critical locus under $\psi\times\bar{\psi}$ coincides with the variety given by the zero locus of the Hermitian critical ideal \eqref{cond} of $(\bu,\bv)$. If the parametrization $\psi$ is generically $d$-to-one then the critical locus in $U\times\overline{U}\subseteq\C^{2m}$ is finite and its cardinality is equal to $d^2\cdot\vHD(X)$.\newline

In the next paragraph we want to give a different intuition of our problem, by rephrasing it in terms of a Euclidean distance problem.\newline

Consider the decomposition $\V=\V_{\R}\oplus i\V_{\R}$, where $\V_{\R}\subseteq\V$ is the subspace invariant under the action of the conjugation, that is the realification of $\V$. Now, note that using coordinates $\bu=\bu^{\Re}+i\bu^{\Im}$ and $\bz=\bx+i\by$ we can write
\begin{equation*}
  q_{\bu}(\bz)=q_{\bu^{\Re}+i\bu^{\Im}}(\bx+i\by)=q_{\bu^{\Re}}(\bx)+q_{\bu^{\Im}}(\by).
\end{equation*}
Split the polynomials $f_1,\ldots,f_s$ defining $X$ into their real and imaginary parts $f_k=f_k^{\Re}+if_k^{\Im}$ such that $f_k^{\Re},f_k^{\Im}\in\R[\bx,\by]$ for $k=1,\ldots,s$ and denote the real algebraic variety
\begin{equation*}
\scalemath{0.9}{\tilde{X}\coloneqq\left\lbrace\left(\frac{\bz+\bw}{2},\frac{\bz-\bw}{2i}\right)\in\V_{\R}^2\otimes\C\mid\bz\in X, \bw\in\overline{X}\right\rbrace=V(f_1^{\Re},f_1^{\Im},\ldots,f_s^{\Re},f_s^{\Im})\subseteq\V_{\R}^2\otimes\C}.
\end{equation*}
The following proposition links our problem to the Euclidean distance problem.

\begin{Prop}\label{rel-com}
Let $X\subseteq\V$ be an algebraic variety. For any $(\bu,\bv)\in\V^2$ the map 
\begin{equation*}
(\bz,\bw)\mapsto\left(\frac{\bz+\bw}{2},\frac{\bz-\bw}{2i}\right)\qquad\text{with inverse}\qquad(\bx,\by)\mapsto(\bx+i\by,\bx-i\by),
\end{equation*}
from $\V^2$ to the complexification $\V_{\R}^2\otimes\C$, sends the points in the zero set of the Hermitian critical ideal \eqref{cond} of $X$ of $(\bu,\bv)$ to the points in the zero set of the critical ideal \eqref{conded} of $\tilde{X}$ of $\left(\frac{\bu+\bv}{2},\frac{\bu-\bv}{2i}\right)$, in particular 
\begin{equation*}
	\vHD(X)=\ED(\tilde{X}).
\end{equation*}
Moreover, if $\bv=\bar{\bu}$ this map sends critical points into critical points, in particular
\begin{equation*}
	\HD(X)=\rED(\tilde{X}).
\end{equation*}
\end{Prop}
\begin{proof}
The first part follows from the linearity of the Wirtinger derivatives and the validity of the Cauchy-Riemann equations for holomorphic function in several variables, by unfolding the definitions of the critical ideals. In fact, the variety $X$ is the zero locus of the $s$ polynomials $f=(f_1,\ldots,f_s)$ and of codimension $c$ iff the variety $\tilde{X}$ is the zero locus of the $2s$ polynomials $(f^{\Re},f^{\Im})$ and of codimension $2c$. In particular, the $(s+1)\times n$ matrices
\begin{equation*}
    \begin{bmatrix}
              \bv-\bw\\
              J_{\bz}(f)
  \end{bmatrix}\qquad\text{and}\qquad\begin{bmatrix}
              \bu-\bz\\
              J_{\bw}(f^{\ast})
  \end{bmatrix}
\end{equation*}
have rank at most $c$ if and only if the $(2s+1)\times 2n$ matrix
\begin{equation*}
    \begin{bmatrix}
              \bu^{\Re}-\bx & \bu^{\Im}-\by\\
              J_{\bx}(f^{\Re}) & J_{\by}(f^{\Re})\\
              J_{\bx}(f^{\Im}) & J_{\by}(f^{\Im})
  \end{bmatrix}
\end{equation*}
has rank at most $2c$.

The second part follows by the definition of the map in the assertion.
\end{proof}

Note that the same map of the last proposition yields a biholomorphism between $X\times\overline{X}$ and $\tilde{X}\subseteq\V_{\R}^2\otimes\C$.

A similar argument applies when $X$ is a parametrized variety.\newline

Most importantly, we can apply all the results known about the Euclidean distance problem to the variety $\tilde{X}$. However, using this approach could often be tedious since dividing the polynomials $f_1,\ldots,f_s$ into their real and imaginary parts requires a non trivial computational effort and the geometry of $\tilde{X}$ could be less clear.\newline

We recall an important result.

\begin{Prop}[Proposition 2.6, \cite{dhost2014}]
Let $X\subseteq\V$ be a variety of codimension $c$ that is cut out by real polynomials $f_1,f_2,\dots,f_s$ of degrees $d_1\geq d_2\geq\ldots\geq d_s$. It holds
\begin{equation*}
    \ED(X)\leq d_1\cdots d_c\sum_{k_1+k_2+\ldots+k_c\leq n-c}(d_1-1)^{k_1}\cdots(d_c-1)^{k_c}
\end{equation*}
and the equality holds for generic varieties.\label{EDvar}
\end{Prop}

As a consequence of the last result, if $X\subseteq\V$ is a variety of codimension $c$ that is cut out by polynomials $f_1,f_2,\dots,f_s$ of degrees $d_1\geq d_2\geq\ldots\geq d_s$, then
\begin{align}
    \vHD(X)&=\ED(\tilde{X})\\
    &\leq (d_1\cdots d_c)^2\sum_{k_1+j_1+\ldots+k_c+j_c\leq 2n-2c}(d_1-1)^{k_1+j_1}\cdots(d_c-1)^{k_c+j_c}\notag\\
    &=(d_1\cdots d_c)^2\sum_{\ell_1+\ldots+\ell_c\leq 2n-2c}(\ell_1+1)\cdots(\ell_c+1)(d_1-1)^{\ell_1}\cdots(d_c-1)^{\ell_c}.\label{firstbound}
\end{align}
Even if we have obtained an upper bound for the value of $\vHD(X)$, this bound is not sharp since, as we discussed above, the variety $\tilde{X}$ is far from being generic. We will sharpen the bound for the value of $\vHD(X)$ of a general hypersurface for several cases in Subsection~\ref{ssec:hyps}. We will see that the true value is much lower.\newline

In the Euclidean distance problem, since all the equations have real coefficients, the conjugate solutions come in pairs. The Hermitian distance problem is in particular a Euclidean distance problem and something similar should happen. This behavior is described in the following proposition.

\begin{Prop}\label{parity}
Let $X$ be an algebraic variety, if $\vHD(X)$ is even (odd) then all the numbers in $\HD(X)$ are even (odd).
\end{Prop}
\begin{proof}
The Hermitian critical ideal of $(\bu,\bar{\bu})$ is invariant under the action of the map $\ast$ and the proposition follows from \cite[Lemma 3.5]{f2024}.
\end{proof}

From this last result it follows that two cases can occur for any solution of the Hermitian critical ideal of $(\bu,\bar{\bu})$. Either it is a conjugated single or, equivalently, a critical point or belongs to an \emph{associated pair} of \cite[Definition 4.6]{f2024}. The latter is a pair of points of the form $(\bz,\bw),(\bar{\bw},\bar{\bz})\in X\times\overline{X}$.

\begin{rmk}\label{realpoints}
A similar argument when $X$ is a real algebraic variety shows that for any solution $(\bz,\bw)$ of the Hermitian critical ideal of $(\bu,\bu)$ we obtain also $(\bw,\bz)$ as another solution. In particular, when $\bu$ is real then for any solution $(\bz,\bw)$ of the Hermitian critical ideal of $(\bu,\bu)$ we also obtain the solutions $(\bw,\bz),(\bar{\bz},\bar{\bw}),(\bar{\bw},\bar{\bz})$. Thus, three cases can occur for any solution. Either it is a real critical point or a complex critical point, which yields another complex critical point, or belongs to one of two associated pairs.
\end{rmk}

We will see that it is easy to construct examples for which the $\HD$ does not contain all the numbers of the same parity between $\min\HD(X)$ and $\max\HD(X)$, see Lemma~\ref{product} below.\newline

\begin{rmk}\label{diltras}
It is not hard to prove that dilatation and translation do not effect the values of the Hermitian Distance degree. 
\end{rmk}

The following result is more interesting to point out.

\begin{Lemma}\label{product}
Let $\W$ denote a $m$-dimensional complex vector space and let $X\subseteq\V$, $Y\subseteq\W$ be algebraic varieties. It holds the equality 
\begin{equation*}
    \vHD(X\times Y)=\vHD(X)\cdot\vHD(Y)
\end{equation*}
and if $\HD(X)=\lbrace a_1,\ldots,a_{d_1}\rbrace$ and $\HD(Y)=\lbrace b_1,\ldots,b_{d_2}\rbrace$ for some $d_1,d_2\in\N$, then 
\begin{equation*}
    \HD(X\times Y)=\lbrace a_1b_1,\ldots,a_1b_{d_2},\ldots,a_{d_1}b_1,\ldots,a_{d_1}b_{d_2}\rbrace.
\end{equation*}
\end{Lemma}

We can use this last result to construct an example of algebraic variety such that the set $\HD$ does not contain all the numbers of the same parity between $\min\HD(X)$ and $\max\HD(X)$.

\begin{ex}\label{exnotcon}
Let $X$ and $Y$ be both the variety of Example~\ref{parab1}. This is the case of a parabola and for it we will compute the values $\vHD(X)=5$ and $\HD(X)=\lbrace 1,3\rbrace$. Using Lemma~\ref{product} there hold the equalities $\vHD(X\times Y)=25$ and $\HD(X\times Y)=\lbrace 1,3,9\rbrace$.\newline
\end{ex}

Let $\Iso(\V)\subseteq\GL(\V)$ be the group of affine linear transformations of $\V$ which preserve the Hermitian inner product. Moreover, for any affine linear transformation $g\in\GL(\V)$ we denote with $\bar{g}$ the affine linear transformation given by the composition $\bar{\phantom{v}}\circ g\circ\bar{\phantom{v}}$.

We provide a particularly useful lemma that characterizes the action of some elements in the group $\Iso(\V)$.

\begin{Lemma}\label{actions}
Let $X\subseteq\V$ be an algebraic variety and $g\in\Iso(\V)$ be an element that leaves $X$ invariant. If $(\bz,\bw)$ is a solution of the Hermitian critical ideal of $(\bu,\bv)$, then $(g\cdot\bz,\bar{g}\cdot\bw)$ is a solution of the Hermitian critical ideal of $(g\cdot\bu,\bar{g}\cdot\bv)$.
\end{Lemma}
\begin{proof}
The assertion follows from the equalities 
\begin{equation*}
  \langle\bu-\bz,\bv-\bw\rangle_{\R}=\langle g\cdot(\bu-\bz),\bar{g}\cdot(\bv-\bw)\rangle_{\R}=\langle g\cdot\bu-g\cdot\bz,\bar{g}\cdot\bv-\bar{g}\cdot\bw\rangle_{\R}.
\end{equation*}
\end{proof}

Using the same proof of the lemma above we can prove the following result.

\begin{Lemma}\label{actionX}
Let $X\subseteq\V$ be an algebraic variety and $g\in\Iso(\V)$. If $(\bz,\bw)$ is a solution of the Hermitian critical ideal of $X$ of $(\bu,\bv)\in\V^2$, then $(g\cdot\bz,\bar{g}\cdot\bw)$ is a solution of the Hermitian critical ideal of $g\cdot X$ of  $(g\cdot\bu,\bar{g}\cdot\bv)\in\V^2$.
\end{Lemma}

We consider now the case of a projective variety.

Let $X\subseteq\V$ be an affine cone and $I_X\subseteq\C[\bz]$ a homogeneous ideal. By a slight abuse of notation, we denote the projective variety in $\Pro\V$ as its affine cone $X$. The $\vHD$ of a projective variety will be the $\vHD$ of the relative affine cone and similarly for the $\HD$.

To take advantage of the homogeneity of the generators of $I_X$, we define
\begin{equation*}
 I_{\Pro X}^{\prime}\coloneqq\left\langle\text{$(c+2)$-minors of $\begin{bmatrix}
        \bv\\
        \bw\\
        J(f)
  \end{bmatrix}$}\right\rangle,
\end{equation*}
and replace the Hermitian critical ideal in expression \eqref{cond} with the following bi-homogeneous ideal
\begin{equation}\label{cond3}
  \left( I_{X}+(I_{X})^{\ast}+I_{\Pro X}^{\prime}+(I_{\Pro X}^{\prime})^{\ast}\right)\colon\left( I_{X_{\sing}}\cdot (I_{X_{\sing}})^{\ast}\cdot \langle\bz,\bw\rangle_{\R}\right)^{\infty}\subseteq\C[\bz,\bw,\bu,\bv].
\end{equation}
Similarly to the general case, we will also refer to the Hermitian critical ideal of $(\bu,\bv)$ when $\bu$ and $\bv$ are fixed.

We denote the projective variety $\tilde{Q}\coloneqq V\left(\langle\bz,\bw\rangle_{\R}\right)\subseteq\V^2$. Observe that the variety $X_{\circ}$ is never contained in $\tilde{Q}$. In particular, any non zero critical point $(\bz,\bar{\bz})$ does not lie in $\tilde{Q}$ since $\langle\bz,\bar{\bz}\rangle_{\R}>0$.\newline

The following lemma concerns the transition between affine cones and projective varieties.

\begin{Lemma}
Let $X\subseteq\V$ be an affine cone and $(\bu,\bv)\in\V^2$ be a point. Let $(\bz,\bw)$ be such that the corresponding point $([\bz],[\bw])$ does not lie in $\tilde{Q}$. The point $([\bz],[\bw])$ lies in the projective variety of the Hermitian critical ideal \eqref{cond3} of $(\bu,\bv)$ if and only if for some unique scalars $\mu_{\bz},\mu_{\bw}\in\C$ the point $(\mu_{\bz}\bz,\mu_{\bw}\bw)$ lies in the affine variety of the Hermitian critical ideal \eqref{cond} of $(\bu,\bv)$. Moreover, when $\bv=\bar{\bu}$ it holds $[\bz]=[\bar{\bw}]$ iff $\mu_{\bz}\bz=\bar{\mu}_{\bw}\bar{\bw}$.
\end{Lemma}
\begin{proof}
Since both Hermitian critical ideals \eqref{cond} and \eqref{cond3} are saturated with respect to $I_{X_{\sing}}\cdot (I_{X_{\sing}})^{\ast}$ and the definitions are symmetric under the map $\ast$, it is sufficient to prove the assertion for $(\bz,\bw)\in  X_{\circ}$ where the Jacobians $J(f)$ at $\bz$ and $J(f^{\ast})$ at $\bw$ have rank $c$.

If $\bv-\mu_{\bw}\bw$ lies in $\Row(J(f))$ at $\mu_{\bz}\bz$, then the subspace $\spn\lbrace\bv,\bw\rbrace+\Row(J(f))$ has dimension at most $c+1$. The rest of the proof of the if condition follows using a symmetric reasoning.

Conversely, suppose that $([\bz],[\bw])$ lies in the variety of the ideal \eqref{cond3}. First assume that $\bw$ lies in $\Row(J(f))$. Thus, $\bw=\sum_{k=1}^s\mu_k\nabla f_k(\bz)$ for some $\mu_k\in\C$ where $k=1,\ldots,s$. Now recall that if $g$ is a homogeneous polynomial in $\C[\bz]$ of degree $d$, then $\langle\bz, \nabla g(\bz)\rangle_{\R}=d\cdot g(\bz)$. Since $f_k(\bz)=0$ for any $k$, we find that $\langle\bz,\nabla f_k(\bz)\rangle_{\R}=0$ for any $k$, which implies that $\langle\bz,\bw\rangle_{\R}=0$. This contradicts our hypothesis, so the matrix
$\begin{bmatrix}
        \bw\\
        J(f)
  \end{bmatrix}$
has rank $c+1$. Thus $\bv-\mu_{\bw}\bw$ lies in $\Row(J(f))$ for a unique $\mu_{\bw}\in\C$. The rest of the proof of the only if condition follows using a symmetric reasoning.

The last assertion can be proved using the argument above by adding the hypothesis of conjugation in each direction.
\end{proof}

From the lemma above and the fact $X_{\circ}\nsubseteq \tilde{Q}$, we directly obtain the following result.

\begin{Coro}
Let $X\subseteq\Pro\V$ be a projective variety and $(\bu,\bv)\in(\Pro\V)^2$ be a generic point. The $\vHD(X)$ is equal to the number of solutions in $(\Pro\V)^2$ of the Hermitian critical ideal $\eqref{cond3}$ of $(\bu,\bv)$ and similarly for $\HD(X)$.
\end{Coro}


\subsection{The Hermitian critical set}\label{ssec:crit}

Let $X\subseteq\V$ be an algebraic variety and $G\subseteq\GL(\V)$ be a Lie subgroup that leaves $X$ invariant. We denote with $\mathfrak{g}=T_{e}G$ the complex Lie algebra of $G$ where $e$ is the identity element. The tangent space to the orbit $G\cdot\bu$ at $\bu\in\V$ is $\bu+\mathfrak{g}\cdot\bu$. Denoting by $G_{\bu}=\lbrace g\in G\mid g\cdot\bu=\bu\rbrace$ the isotropy group of $\bu$, we have $\dim\mathfrak{g}\cdot\bu=\dim\mathfrak{g}-\dim G_{\bu}$.\newline

In \cite{o2022}, assuming that the elements of $G$ preserve a positive symmetric bilinear form $\hat{q}$ on $\V$ for any $\bu\in\V$ is defined the \emph{critical space} that is the vector subspace $\lbrace\bv\in\V\mid\hat{q}(\bv,g\cdot\bu)=0\quad\forall g\in\mathfrak{g}\rbrace$. The key property of this space is that it contains the critical points of the distance function from $X$ induced by $\hat{q}$. We try to adapt similar arguments in the Hermitian case.\newline

We recall that $q$ indicates a Hermitian form on $\V$ and let $\U(\V)\subseteq\GL(\V)$ be the group of linear transformations of $\V$ which preserve $q$.

\begin{Lemma}\label{lemcrit}
Let $G\subseteq\U(\V)$, if $g\in\mathfrak{g}$ then $q(g\cdot \bu,\bv)=-q(\bu,g\cdot\bv)$ for any $\bu,\bv\in\V$ and in particular $q(g\cdot\bu,\bu)^{\Re}=0$.
\end{Lemma}
\begin{proof}
Let $\psi(t)\colon[-1,1]\to G$ be a path such that $\psi(0)=e$ and $\psi^{\prime}(0)=g$. Taking the
derivative at $t=0$ of the constant function $q(\psi(t)\cdot\bu,\psi(t)\cdot\bv)=q(\bu,\bv)$ the claim follows.

The second assertion follows from the chain of equalities
\begin{equation*}
    q(g\cdot\bu,\bu)=-q(\bu,g\cdot\bu)=-\overline{q(g\cdot\bu,\bu)}.
\end{equation*}
\end{proof}

If we consider a symmetric bilinear form $\hat{q}$ and a group $G$ that preserves it, the same proof of the last lemma shows $\hat{q}(g\cdot\bu,\bu)=0$. This is the main difference between the symmetric and the Hermitian cases, resulting in a more complicated definition in the latter since the set we focus on is not an affine subspace of $\V$.

\begin{Def}
Let $\bu\in\V$, we define the \emph{Hermitian critical set} of $\bu$ as the subset
\begin{equation*}
    H_{\bu}\coloneqq\lbrace\bv\in\V\mid q(\bv,g\cdot(\bv-\bu))=0\quad\forall g\in\mathfrak{g}\rbrace.
\end{equation*}
\end{Def}

\begin{rmk}
We can notice from Lemma~\ref{lemcrit} that for $\bv\in H_{\bu}$ there hold
\begin{equation*}
	0=q(\bv,g\cdot(\bv-\bu))^{\Re}=q(\bv,g\cdot\bv)^{\Re}-q(\bv,g\cdot\bu)^{\Re}=-q(\bv,g\cdot\bu)^{\Re}.
\end{equation*}
In particular, the condition $q(\bv,g\cdot\bu)=\rho e^{\pm\frac{\pi}{2}i}$ for some $0\leq\rho\in\R$, or in other words that the Kasner’s pseudo-angle of $\bv$ and $g\cdot\bu$ is equal to $\pm\pi/2$, defines a linear subspace of $\V$ that contains the Hermitian critical set.
\end{rmk}

\begin{Theo}\label{teo1crit}
Let $X\subseteq\V$ be an algebraic variety and $G\subseteq\U(\V)$ be a subgroup that leaves $X$ invariant. The regular critical points of $q_{\bu}$ lie in $H_{\bu}$.
\end{Theo}
\begin{proof}
Let $\bz\in X$ be a critical point, from the fact that $X$ is $G$ invariant follows $\mathfrak{g}\cdot\bz\subseteq T_{\bz}X$. Thus the result follows from the equalities
\begin{equation*}
    0=q(g\cdot\bz,\bz-\bu)=-q(\bz,g\cdot(\bz-\bu))
\end{equation*}
where on the left we used the fact that $\bz$ is critical and on the right Lemma~\ref{lemcrit}.
\end{proof}

\begin{Theo}
Let $X\subseteq\V$ be an algebraic variety, $G\subseteq\U(\V)$ be a subgroup that leaves $X$ invariant and $\bz\in H_{\bu}\cap X$. There hold the following
\begin{enumerate}[\normalfont i)]
    \item if the orbit $G\cdot\bz$ is dense in $X$ then $\bz$ is a critical point of $q_{\bu}$,
    \item if $X$ is an affine cone and the orbit $G\cdot[\bz]$ is dense in the projective variety $X\subseteq\Pro\V$ then there exists a unique $\lambda\in\C$ such that $\lambda\bz$ is a critical point of $q_{\bu}$.
\end{enumerate}
\end{Theo}
\begin{proof}
By assumption we have the equality $\mathfrak{g}\cdot\bz=T_{\bz}X$ and since the steps of the proof of Theorem~\ref{teo1crit} are invertible this proves i).

The assumption of ii) implies that $\mathfrak{g}\cdot\bz+\spn\lbrace\bz\rbrace=T_{\bz}X$. Since $q(\bz,\bz)$ is non zero for $\lambda=q(\bu,\bz)/q(\bz,\bz)$ it holds $q(\lambda\bz,\lambda\bz-\bu)=0$, so that we get orthogonality on $\spn\lbrace\bz\rbrace\subseteq T_{\bz}X$. Orthogonality to $\mathfrak{g}\cdot\bz$ then follows by applying the same argument of i) and replacing $\bz$ by $\lambda\bz$.
\end{proof}


\subsection{vHDdeg of hypersurfaces}\label{ssec:hyps}

We discuss the value of $\vHD(X)$ when $X$ is an algebraic variety of codimension one. The results obtained in this section are easily achieved thanks to the clear symmetries of the problem given by the complex approach.

Let $X=V(f)\subseteq\V$ be a hypersurface of degree $d$. By means of equation \eqref{firstbound}, we already know 
\begin{align}\label{firstboundhyp}
    \vHD(X)&\leq d^2\sum_{k=0}^{2n-2}(k+1)(d-1)^k\\
    &=\frac{d^2\left((2n-2)(d-1)^{2n}-(2n+1-d)(d-1)^{2n-1}+1\right)}{(d-2)^2}.\notag
\end{align}
In particular, to obtain the second formula we assume $d>2$. 

We recall that if $X$ is a real hypersurface, or in other words $f$ is a real polynomial, from Proposition~\ref{EDvar} we get $\ED(X)\leq d\sum_{k=0}^{n-1}(d-1)^k$ and the equality holds for generic hypersurface.

\begin{Prop}\label{deghyps}
Let $X=V(f)\subseteq\V$ be a generic hypersurface of degree $d$, it holds 
\begin{equation}\label{secondboundhyp}
    \vHD(X)\leq d^2\sum_{k=0}^{n-1}\binom{n-1}{k}^2(d-1)^{2k}.
\end{equation}
\end{Prop}
\begin{proof}
The codimension of $X$ is $c=1$. The number of solutions of the Hermitian critical ideal of $X$ is bounded using the Bernstein–Kushnirenko theorem by the value
\begin{equation*}
   \MV_{2n}(\triangle_f,\triangle_{f^{\ast}},\triangle_{g_{1}},\ldots,\triangle_{g_{n-1}},\triangle_{g^{\ast}_{1}},\ldots,\triangle_{g^{\ast}_{n-1}}) 
\end{equation*}
where the first two arguments are the convex hulls 
\begin{alignat*}{2}
    &\triangle_f\coloneqq\conv\lbrace(\bs,\bo)\in(\N^n)^2\mid\|\bs\|_{1}\leq d\rbrace,\qquad &&\triangle_{f^{\ast}}\coloneqq\conv\lbrace(\bo,\bt)\in(\N^n)^2\mid\|\bt\|_{1}\leq d\rbrace,
\end{alignat*}
while the remaining are in the families of convex hulls
\begin{align*}
    &\triangle_{g_{\ell,r}}\coloneqq\conv\lbrace(\bs,\bt)\in(\N^n)^2\mid\|\bs\|_{1}\leq d-1,\ \|\bt\|_{1}\leq 1,\ t_{\xi}=0\text{ if }\xi\neq\ell,r\rbrace,\\
    &\triangle_{g_{\ell,r}^{\ast}}\coloneqq\conv\lbrace(\bs,\bt)\in(\N^n)^2\mid\|\bs\|_{1}\leq 1,\ \|\bt\|_{1}\leq d-1,\ s_{\xi}=0\text{ if }\xi\neq\ell,r\rbrace
\end{align*}
for $\ell<r\in[n]$ respectively. The first two come from the equations $f=f^{\ast}=0$ while the other $2\binom{n}{2}$ come from the equations given by the order $2$ minors of the matrices involved in the Hermitian critical ideal and only $2n-2$ of them are required since our matrices are $2\times n$. Denote the convex hulls
\begin{equation*}
   \triangle_{1,0}\coloneqq\conv(\triangle\times\lbrace\bo\rbrace)\subseteq\R^{2n}\qquad\text{and}\qquad\triangle_{0,1}\coloneqq\conv(\lbrace\bo\rbrace\times\triangle)\subseteq\R^{2n},
\end{equation*}
where $\triangle\coloneqq\lbrace\bs\in\N^n\mid\|\bs\|_{1}\leq 1\rbrace\subseteq\R^{n}$, so that 
\begin{equation*}
    \triangle_{f}=d\triangle_{1,0}\subseteq\R^{2n}\quad\text{and}\quad\triangle_{f^{\ast}}=d\triangle_{0,1}\subseteq\R^{2n}.
\end{equation*}
Moreover, there hold the inclusions
\begin{equation*}
    \triangle_{g_{\ell,r}}\subseteq\triangle_{\mathcal{Z}}\coloneqq(d-1)\triangle_{1,0}+\triangle_{0,1}\subseteq\R^{2n}\quad\text{and}\quad\triangle_{g^{\ast}_{\ell,r}}\subseteq\triangle_{\mathcal{W}}\coloneqq\triangle_{1,0}+(d-1)\triangle_{0,1}\subseteq\R^{2n}.
\end{equation*}
It will be clear by the rest of the proof that the mixed volume we are computing is equal to
\begin{equation*}
   \MV_{2n}(\triangle_f,\triangle_{f^{\ast}},\underbrace{\triangle_{\mathcal{Z}},\ldots,\triangle_{\mathcal{Z}}}_{n-1},\underbrace{\triangle_{\mathcal{W}},\ldots,\triangle_{\mathcal{W}}}_{n-1}).
\end{equation*}
We make this change in order to simplify the notation. Since the mixed value is symmetric multilinear, we rewrite this last value as
\begin{equation*}
    \scalemath{0.9}{d^2\sum_{k=0}^{n-1}\binom{n-1}{k}\sum_{j=0}^{n-1}\binom{n-1}{j}(d-1)^{n-1+k-j}\MV_{2n}(\triangle_{1,0},\triangle_{0,1}\underbrace{\triangle_{1,0},\ldots,\triangle_{1,0}}_{k+j},\underbrace{\triangle_{0,1},\ldots,\triangle_{0,1}}_{2n-2-k-j})}.
\end{equation*}
Now, \cite[Lemma 4.5]{e1996} implies that the mixed volumes of the terms with $k+j\neq 2n-2-k-j$, or equivalently $k\neq n-1-j$, vanish, thus the formula simplifies to
\begin{equation*}
    d^2\sum_{k=0}^{n-1}\binom{n-1}{k}^2(d-1)^{2k}\MV_{2n}(\underbrace{\triangle_{1,0},\ldots,\triangle_{1,0}}_{n},\underbrace{\triangle_{0,1},\ldots,\triangle_{0,1}}_{n}).
\end{equation*}
Applying the same lemma for the remaining terms we get
\begin{align*}
    \MV_{2n}(\underbrace{\triangle_{1,0},\ldots,\triangle_{1,0}}_{n},\underbrace{\triangle_{0,1},\ldots,\triangle_{0,1}}_{n})&=\MV_{n}(\triangle,\ldots,\triangle)\MV_{n}(\triangle,\ldots,\triangle)=1.
\end{align*}
If we do not change the convex bodies in the formula, the steps will be the same until the last equality. In this case, on the left hand side, instead of the standard simplex $\triangle\subseteq\R^n$, we could have instances of standard simplexes in pairwise distinct coordinate subspaces $\R^{n-1}\subseteq\R^n$. We can replace those convex bodies with $\triangle$ by means of \cite[Theorem 2]{c2019} and the claim follows.
\end{proof}

The symmetric degrees in $\bz$ and $\bw$ is not the only property shared by the polynomials of the system given by the Hermitian critical ideal. In fact, the Hermitian critical ideal is invariant under the action of the map $\ast$, thus the true bound of $\vHD$ could be lower, see Remark~\ref{bounds}. Nonetheless, we already said that $\max\HD$ could be lower than $\vHD$, see the examples in Subsection~\ref{ssec:conics}.

\begin{rmk}\label{bounds}
The bound of equation \eqref{secondboundhyp} is not sharp in general. In fact, the bound of equation \eqref{firstboundhyp} is still better for $n$ tending to infinity. For example, assume $X$ to be a hypersurface of degree $d=2$. In this case equation \eqref{secondboundhyp} yields the value 
\begin{equation*}
    4\sum_{k=0}^{n-1}\binom{n-1}{k}^2=4\binom{2n-2}{n-1}\geq 4\cdot 2^{n-1},
\end{equation*}
where we used the Chu-Vandermonde identity and a well known inequality on the binomial coefficient $\binom{2n-2}{n-1}$ to get a lower bound on the right. On the other hand, equation \eqref{firstboundhyp} yields the value $4\sum_{k=0}^{2n-2}(k+1)=4\cdot n(2n-1)$. In particular, the value of \eqref{firstboundhyp} is is surely lower than the value of \eqref{secondboundhyp} for $n\geq 8$ when $d=2$.
\end{rmk}

The last remark motivates the following result.

\begin{Coro}
Let $X=V(f)\subseteq\V$ be a generic hypersurface of degree $d$, it holds
\begin{equation}\label{finalboundhyp}
    \vHD(X)\leq d^2\min\left\lbrace\sum_{k=0}^{n-1}\binom{n-1}{k}^2(d-1)^{2k},\ \sum_{k=0}^{2n-2}(k+1)(d-1)^k\right\rbrace.
\end{equation}
\end{Coro}

We compute some values of the upper bound of the result above in Table~\ref{tabhyps}. Clearly, for the case $d=1$ it holds $\vHD=1$ for any $n$.

\begin{table}[h!]
\caption{Values of equation \eqref{finalboundhyp} for different $d$ and $n$. The underlined values are the ones where the minimum is reached by the formula on the right. The values with a hat are the ones we checked to be sharp.}\label{tabhyps}
  \begin{tabular}{c|cccccccc}
    $d$ & $n=1$ & $n=2$ & $n=3$ & $n=4$ & $n=5$ & $n=6$ & $n=7$ & $n=8$\\
    \hline
    $2$ & $\hat{4}$ & $\hat{8}$ & $24$ & $80$ & $\underline{180}$ & $\underline{264}$ & $\underline{364}$ & $\underline{480}$\\
    $3$ & $\hat{9}$ & $\hat{45}$ & $297$ & $2205$ & $17289$ & $139725$ & $\underline{884745}$ & $\underline{4128777}$\\
    $4$ & $\hat{16}$ & $\hat{160}$ & $1888$ & $24640$ & $340576$ & $4868800$ & $71097280$ & $1053289600$\\
    $5$ & $\hat{25}$ & $\hat{425}$ & $8025$ & $163625$ & $3513625$ & $78064425$ & $1774203225$ & $40958848425$\\
  \end{tabular}
\end{table}

In the following example we present the $\vHD$ of a curve of degree $3$ and of the Fermat cubic.

\begin{ex}\label{exfcubic}
The curve 
\begin{equation*}
    X=V(z_1^3+z_2^3-1)\subseteq\C^2
\end{equation*}
satisfies $\vHD(X)=45=3^2(2\cdot 2+1)$ as predicted by Proposition~\ref{deghyps}, this can be seen by computing the degree of the Hermitian critical ideal.\newline

(Fermat cubic). The Fermat cubic 
\begin{equation*}
    X=V(z_1^3+z_2^3+z_3^3-1)\subseteq\C^3
\end{equation*}
satisfies $\vHD(X)=189=3^3(2\cdot 3+1)$, this can be seen by calculating the degree of the Hermitian critical ideal. This value is lower than $3^3\cdot 11=297$ predicted by Proposition~\ref{deghyps}.
\end{ex}

Related to these last computations we state the following.

\begin{Conj}
Let $2\leq n\in\N$ and
\begin{equation*}
    X=V(z_1^3+z_2^3+\ldots+z_n^3-1)\subseteq\C^n,
\end{equation*}
it holds $\vHD(X)=3^n(2n+1)$.
\end{Conj}


\subsection{vHDdeg of parametrized varieties}\label{ssec:param}

We discuss the value of $\vHD(X)$ when $X$ admits a polynomial parametrization.

As in Subsection~\ref{ssec:hyps}, the results obtained in this section are easily achieved thanks to the clear symmetries of the problem given by the complex approach.

Let $X\subseteq\C^n$ be an algebraic variety parametrized by polynomials where $m$ of them are generic of degree $d$ and the remaining are generic of degree less or equal to $d$. The B\'ezout's theorem applied to the system of polynomial equations $\nabla_{\bz,\bw} D_{\bu,\bv}(\bz,\bw)=\bo$ yields the bound 
\begin{equation}\label{firstboundparam}
    \vHD(X)\leq(2d-1)^{2m}.
\end{equation}

On the other hand, if $X$ is a real algebraic variety then a similar argument shows that it holds $\ED(X)\leq(2d-1)^m$ and the equality holds for generic parametrization.

\begin{Prop}\label{degparam}
Let $X$ be a generic algebraic variety parametrized by $n$ polynomials of degree $d$ in $m$ variables, it holds
\begin{equation}\label{secondboundparam}
    \vHD(X)\leq\sum_{k=0}^{m}\binom{m}{k}^{2}d^{2(m-k)}(d-1)^{2k}.
\end{equation}
\end{Prop}
\begin{proof}
Let $m$ polynomials be generic of degree $d$ and the remaining be generic of degree less than or equal to $d$. The number of solutions of the system $\nabla_{\bz,\bw}D_{\bu,\bv}(\bz,\bw)=\bo$ is bounded using the Bernstein–Kushnirenko theorem by the value 
\begin{equation*}
    \MV_{2m}(\underbrace{\triangle_{\mathcal{Z}},\ldots\triangle_{\mathcal{Z}}}_m,\underbrace{\triangle_{\mathcal{W}},\ldots,\triangle_{\mathcal{W}}}_m).
\end{equation*}
The first $m$ arguments of this mixed volume come from the $m$ equations $\nabla_{\bz}D_{\bu,\bv}(\bz,\bw)=\bo$, and are all the same convex hull $\triangle_{\mathcal{Z}}\coloneqq\conv\mathcal{Z}\subseteq\R^{2m}$, and the second $m$ arguments come from the $m$ equations $\nabla_{\bw}D_{\bu,\bv}(\bz,\bw)=\bo$, and are all the same convex hull $\triangle_{\mathcal{W}}\coloneqq\conv\mathcal{W}\subseteq\R^{2m}$, where we used the two sets $\mathcal{Z}\coloneqq\lbrace (\bs,\bt)\in(\N^m)^2\mid \|\bs\|_{1}\leq d-1, \|\bt\|_{1}\leq d\rbrace$ and $\mathcal{W}\coloneqq\lbrace (\bt,\bs)\mid (\bs,\bt)\in\mathcal{Z}\rbrace$ respectively. Denote the convex hulls
\begin{equation*}
   \triangle_{1,0}\coloneqq\conv(\triangle\times\lbrace\bo\rbrace)\subseteq\R^{2m}\qquad\text{and}\qquad\triangle_{0,1}\coloneqq\conv(\lbrace\bo\rbrace\times\triangle)\subseteq\R^{2m},
\end{equation*}
where $\triangle\coloneqq\lbrace\bs\in\N^m\mid\|\bs\|_{1}\leq 1\rbrace\subseteq\R^{m}$, so that 
\begin{equation*}
    \triangle_{\mathcal{Z}}=(d-1)\triangle_{1,0}+d\triangle_{0,1}\quad\text{and}\quad\triangle_{\mathcal{W}}=d\triangle_{1,0}+(d-1)\triangle_{0,1}.
\end{equation*}
Thus, using the fact that the mixed volume is symmetric multilinear, we can rewrite the value $\MV_{2m}(\triangle_{\mathcal{Z}},\ldots\triangle_{\mathcal{Z}},\triangle_{\mathcal{W}},\ldots,\triangle_{\mathcal{W}})$ as
\begin{equation*}
    \sum_{k=0}^m\binom{m}{k}\sum_{j=0}^m\binom{m}{j}(d-1)^{m+k-j}d^{m+j-k}\MV_{2m}(\underbrace{\triangle_{1,0},\ldots,\triangle_{1,0}}_{k+j},\underbrace{\triangle_{0,1},\ldots,\triangle_{0,1}}_{2m-k-j}).
\end{equation*}
Now, \cite[Lemma 4.5]{e1996} implies that the mixed volumes of the terms with $k+j\neq 2m-k-j$, or equivalently $k\neq m-j$, vanish, thus the formula simplifies to
\begin{equation*}
    \sum_{k=0}^m\binom{m}{k}^2(d-1)^{2k}d^{2(m-k)}\MV_{2m}(\underbrace{\triangle_{1,0},\ldots,\triangle_{1,0}}_{m},\underbrace{\triangle_{0,1},\ldots,\triangle_{0,1}}_{m}).
\end{equation*}
Applying the same lemma for the remaining terms we get
\begin{align*}
    \MV_{2m}(\underbrace{\triangle_{1,0},\ldots,\triangle_{1,0}}_{m},\underbrace{\triangle_{0,1},\ldots,\triangle_{0,1}}_{m})&=\MV_{m}(\triangle,\ldots,\triangle)\MV_{m}(\triangle,\ldots,\triangle)=1
\end{align*}
and the claim follows.
\end{proof}

\begin{rmk}\label{bounds2}
The bound of equation \eqref{secondboundparam} is always sharper than the bound of equation \eqref{firstboundparam}. In fact, we bound the value of equation \eqref{secondboundparam} by
\begin{align*}
    \sum_{k=0}^{m}\binom{m}{k}^{2}d^{2(m-k)}(d-1)^{2k}&\leq\left(\sum_{k=0}^{m}\binom{m}{k}d^{m-k}(d-1)^{k}\right)^2\\
    &=\left((d+(d-1))^{m}\right)^{2}=(2d-1)^{2m},
\end{align*}
where the last formula of the chain above is exactly the formula of equation \eqref{firstboundparam}.
\end{rmk}

We compute some values of equation \eqref{secondboundparam} in Table~\ref{tabparam}. Clearly, for the case $d=1$ it holds $\vHD=1$ for any $m$. 

\begin{table}[h!]
\caption{Values of equation \eqref{secondboundparam} for different $d$ and $m$. The values with a hat are the ones we checked to be sharp.}\label{tabparam}
  \begin{tabular}{cccccc}
    & \multicolumn{1}{c}{$m=1$} & & \multicolumn{1}{c}{$m=2$} & & \multicolumn{1}{c}{$m=3$}\\
    \cmidrule{2-2}
    \cmidrule{4-4}
    \cmidrule{6-6}
    $\scalemath{0.9}{d}$ & $\scalemath{0.9}{2d^2-2d+1}$ & &  $\scalemath{0.9}{6d^4-12d^3+10d^2-4d+1}$ & & $\scalemath{0.9}{20d^6-60d^5+78d^4-56d^3+24d^2-6d+1}$ \\
    \midrule
    $2$ & $\hat{5}$ & & $33$ & & $245$\\
    $3$ & $\hat{13}$ & & $241$ & & $5005$\\
    $4$ & $\hat{25}$ & & $913$ & & $37225$\\
    $5$ & $\hat{41}$ & & $2481$ & & $167321$\\
  \end{tabular}
\end{table}

One more time, the symmetric degrees in $\bz$ and $\bw$ is not the only property shared by the polynomials of the system $\nabla_{\bz,\bw}D_{\bu,\bv}(\bw,\bz)=\bo$. In fact, this system is invariant under the action of the map $\ast$ and the true value of $\vHD(X)$ could be lower. Moreover, as already said, the number $\max\HD(X)$ could be lower than $\vHD$, see Example~\ref{parab1}.


\subsection{HDdeg of conics}\label{ssec:conics}

In this subsection we focus on the case of conics and provide some results for the $\HD$ of projective and affine conics.\newline

We start by considering projective conics. Thus, we set 
\begin{equation*}
    X=V(az_1^2+2\cdot bz_1z_2+cz_2^2+2\cdot dz_1z_3+2\cdot ez_2z_3+fz_3^2)\subseteq\Pro^2
\end{equation*}
where $a,b,c,d,e,f\in\C$. In particular, $X$ is given by the zero set of the symmetric bilinear form
\begin{equation*}
    \begin{bmatrix}
        z_1 & z_2 & z_3
    \end{bmatrix}
    \begin{bmatrix}
        a & b & d\\
        b & c & e\\
        d & e & f
    \end{bmatrix}
    \begin{bmatrix}
        z_1\\
        z_2\\
        z_3
    \end{bmatrix}=\bz^TA\bz
\end{equation*}
that is uniquely determined by the complex symmetric matrix $A$.

Let us recall a well known fact on complex symmetric matrices that can be found in \cite[Corollary 4.4.4 (c)]{hj2013}.

\begin{Prop}\label{unmat}
Let $A$ be a complex symmetric matrix. There exists a unitary matrix $U$ such that the product $UAU^T$ is a non negative diagonal matrix whose diagonal entries are the singular values of $A$, in any desired order.
\end{Prop}

We will use the result above to limit the number of cases we have to address.

\begin{Theo}\label{projcon}
Let $X$ be a projective conic and $A$ be its representative complex symmetric matrix.

If $A$ is singular then
\begin{enumerate}[\normalfont i)]
  \item $\vHD(X)=1$ and $\HD(X)=\lbrace 1\rbrace$ if $\rk(A)=1$,
  \item $\vHD(X)=2$ and $\HD(X)=\lbrace 2\rbrace$ if $\rk(A)=2$.
\end{enumerate}

If $A$ is non singular then
\begin{enumerate}[\normalfont i)]
  \item $\vHD(X)=2$ and $\HD(X)=\lbrace 2\rbrace$ if $A$ has only one singular value,
  \item $\vHD(X)=6$ if $A$ has two different singular values,
  \item $\vHD(X)=8$ if $A$ has three different singular values.
\end{enumerate}
\end{Theo}
\begin{proof}
By Lemma~\ref{actionX}, since the Hermitian inner product is invariant under the action of unitary matrices on the complex vector space, we can study of the Hermitian distance problem for $X$ on the variety $U\cdot X$ where $U$ is a $3\times 3$ unitary matrix. We combine this fact with Proposition~\ref{unmat} and choose a $3\times 3$ unitary matrix $U$ such that relabeling the parameters and applying a change of variables, the quadratic form representing $U\cdot X$ is
\begin{equation*}
    \bz^TA\bz=\bz^TU^T\bar{U}AU^HU\bz=(U\bz)^T\begin{bmatrix}
        a & 0 & 0\\
        0 & c & 0\\
        0 & 0 & f
    \end{bmatrix}(U\bz)
\end{equation*}
where $0\leq f\leq c\leq a\in\R$ are the singular values of the starting matrix $A$. Moreover, assuming we are not in the trivial case of all vanishing singular values, without loss of generality, we use Remark~\ref{diltras} and by applying the symmetric matrix $\sqrt{a}I_3$, we replace $X$ with the variety $\sqrt{a}U\cdot X$, so that we can also set $a=1$. 

Suppose $A$ is singular. If two singular values are equal to zero, say $c=f=0$, then $X=V(z_1^2)$ is an hyperplane and point i) follows from Proposition~\ref{affsub}. If only one singular value is equal to zero, say $f=0$, then $z_1^2+cz_2^2=(z_1+i\sqrt{c}z_2)(z_1-i\sqrt{c}z_2)$ so that $X$ is the union of two lines and point ii) follows similarly to point i). 

Suppose now that $A$ is non singular. The values of $\vHD(X)$ for all points i), ii) and iii) follow by computing the degree of the Hermitian critical ideal of $X$, for example using Macaulay2. The only value in $\HD(X)$ for point i) follows from the properties of the $\HD$.
\end{proof}

We now consider affine conics. Thus, we set
\begin{equation*}
    X=V(az_1^2+2\cdot bz_1z_2+cz_2^2+2\cdot dz_1+2\cdot ez_2+f)\subseteq\C^2
\end{equation*}
where $a,b,c,d,e,f\in\C$. Similarly to the projective case, the variety $X$ is uniquely determined by the same complex symmetric matrix $A$ studied above.

\begin{Theo}\label{affcon}
Let $X$ be an affine conic, $A$ be its representative complex symmetric matrix and $A_1$ be the $2\times 2$ leading principal submatrix of $A$.

If $A$ is singular
\begin{enumerate}[\normalfont i)]
    \item $\vHD(X)=1$ and $\HD(X)=\lbrace 1\rbrace$ if $\rk(A)=1$,
    \item $\vHD(X)=2$ and $\HD(X)=\lbrace 2\rbrace$ if $\rk(A)=2$.
\end{enumerate}

If $A$ is non singular then
\begin{enumerate}[\normalfont i)]
    \item $\vHD(X)=5$ and $\HD(X)=\lbrace 1,3\rbrace$ if $A_1$ is singular,
    \item $\vHD(X)=6$ and $\HD(X)=\lbrace 2,4\rbrace$ if $A_1$ is non singular and has two equal singular values,
    \item $\vHD(X)=8$ and $\lbrace 2,4\rbrace\subseteq\HD(X)\subseteq\lbrace 2,4,6\rbrace$ if $A_1$ is non singular and has two different singular values.
\end{enumerate}
\end{Theo}
\begin{proof}
If $A$ is singular then
\begin{equation*}
    az_1^2+2bz_1z_2+cz_2^2+2dz_1+2ez_2+f=(\alpha_1z_1+\beta_1z_2+\gamma_1)(\alpha_2z_1+\beta_2z_2+\gamma_2)
\end{equation*}
for some $\alpha_k,\beta_k,\gamma_k\in\C$ for $k=1,2$ so that $X$ is the union of two lines. Moreover, the two lines coincide exactly when $\rk(A)=1$. Points i) and ii) thus follow from Proposition~\ref{affsub}.

In general, similarly to the start of the proof of Theorem~\ref{projcon}, by relabeling the parameters and applying a change of variables, we replace $X$ with the variety $U\cdot X$, where
\begin{equation*}
    U=\begin{bmatrix}
        U_1 & \bo\\
        \bo & 1
    \end{bmatrix}
\end{equation*}
with $U_1$ a $2\times 2$ unitary matrix such that
\begin{equation*}
    \begin{bmatrix}
        z_1 & z_2 & 1
    \end{bmatrix}A\begin{bmatrix}
        z_1\\
        z_2\\
        1
    \end{bmatrix}=\begin{bmatrix}
        \bz\\
        1
    \end{bmatrix}^TU^T\bar{U}AU^HU\begin{bmatrix}
        \bz\\
        1
    \end{bmatrix}=\left(U\begin{bmatrix}
        \bz\\
        1
    \end{bmatrix}\right)^T\begin{bmatrix}
        a & 0 & d\\
        0 & c & e\\
        d & e & f
    \end{bmatrix}\left(U\begin{bmatrix}
        \bz\\
        1
    \end{bmatrix}\right)
\end{equation*}
where $0\leq c\leq a\in\R$ are the singular value of the starting submatrix $A_1$ and $d,e,f\in\C$. With these choices it holds
\begin{equation*}
    \det(A)=\det\begin{bmatrix}
        a & 0 & d\\
        0 & c & e\\
        d & e & f
    \end{bmatrix}=acf-ae^2-cd^2,
\end{equation*}
thus if $A$ is non singular at least one parameter among $a$ and $c$ must be non zero. Without loss of generality say $a\neq 0$. Using Remark~\ref{diltras}, we can translate by the vector $(d/a,0)\in\C^2$ to set $d=0$ and scale by the factor $\sqrt{a}\in\C$, that is the same of applying the symmetric matrix
\begin{equation*}
    \begin{bmatrix}
        \sqrt{a}I_2 & \bo\\
        \bo & 1
    \end{bmatrix}\quad\text{to the vector}\quad\begin{bmatrix}
        \bz\\
        1
    \end{bmatrix},
\end{equation*}
to set $a=1$. In the end, we are left with
\begin{equation*}
    A=\begin{bmatrix}
        1 & 0 & 0\\
        0 & c & e\\
        0 & e & f
    \end{bmatrix}
\end{equation*}
where $0\leq c\leq1$ and $e,f\in\C$. Since translations and scalings do not affect the ratio between the singular values of $A_1$ and, in particular, do not affect the vanishing of $\det(A_1)=c$, we can consider two mutually exclusive cases depending on the parameter $c$. Case (1) will be $c=0$ and case (2) will be $c\neq 0$.

In case (1) we are in the hypothesis $\det(A)=-e^2\neq 0$. Thus, we translate by the vector $\bb=(0,f/2e)\in\C^2$ to set $f=0$ and scale by the factor $-1/2e\in\C$, that is the same of applying the symmetric matrix
\begin{equation*}
    \begin{bmatrix}
        -\frac{1}{2e}I_2 & \bo\\
        \bo & 1
    \end{bmatrix}\quad\text{to the vector}\quad\begin{bmatrix}
        \bz\\
        1
    \end{bmatrix},
\end{equation*}
so that we can multiply the equation obtained by $1/4e^2\in\C$ to set $e=-1/2$ and replace our variety with $X=V(z_1^2-z_2)$. In the end, point i) follows from Example~\ref{parab1}.

In case (2) we apply the translation by the vector $(0,e/c)\in\C^2$ to set $e=0$. Moreover, since we are in the hypothesis $\det(A)=cf\neq 0$, it holds $f\neq 0$. Thus, we scale by the factor $i\sqrt{c/f}\in\C$, that is the same of applying the symmetric matrix
\begin{equation*}
    \begin{bmatrix}
        i\sqrt{c/f}I_2 & \bo\\
        \bo & 1
    \end{bmatrix}\quad\text{to the vector}\quad\begin{bmatrix}
        \bz\\
        1
    \end{bmatrix},
\end{equation*}
so that we can multiply the equation obtained by $-c/f\in\C$ to set $f=-c$ and replace our variety with $X=V(z_1^2+cz_2^2-c)$ where $0<c\leq 1$. If $c=1$ that means the starting submatrix $A_1$ has two equal singular values and point ii) follows from Example~\ref{circle1}. If $c<1$ that means the starting submatrix $A_1$ has two different singular values and point iii) follows from Example~\ref{ellipse1}.
\end{proof}

The case iii) of $A$ non singular remains still incomplete.


\subsection{HD correspondence and duality}\label{ssec:corrdual}

Let $X\subseteq\V$ be an algebraic variety, we define the \emph{virtual Hermitian Distance correspondence} (vHD correspondence) of $X$ denoted $\vH_{X}$ as the subvariety of $\V^2\times\V^2$ defined by the Hermitian critical ideal \eqref{cond}. In particular, it holds the containment $\vH_{X}\subseteq X\times\overline{X}\times\V^2$. The following theorem is the generalization of Lemma~\ref{lemmastart}. 

\begin{Theo}\label{corr}
Let $X\subseteq\V$ be an algebraic variety of codimension $c$. The vHD correspondence $\vH_{X}$ is a $2n$ dimensional irreducible variety. The projection $\pi_{X\times\overline{X}}$ on the first $2$ components is a vector bundle of rank $2c$. The projection $\pi_{\V^2}$ on the second $2$ components has generic fibers of cardinality equal to $\vHD(X)$. Moreover, $\min\HD(X)$ is positive.
\end{Theo}
\begin{proof}
The first part follows as \cite[Theorem 4.1]{dhost2014}. The only remaining part is the fact that the projection into the third component is a dominant map.

Firstly observe that the double diagonal $\Delta(X)\coloneqq\lbrace(\bz,\bar{\bz},\bz,\bar{\bz})\in\V^2\times\V^2\mid \bz\in X\rbrace$ is contained in $\vH_{X}$. Fix a point $\bz\in X_{\reg}$ so that $\vH_{X}$ is smooth at the point $(\bz,\bar{\bz},\bz,\bar{\bz})$. The tangent space $T_{(\bz,\bar{\bz},\bz,\bar{\bz})}\vH_{X}$ contains both the tangent space $T_{(\bz,\bar{\bz},\bz,\bar{\bz})}\Delta(X)=\Delta(T_{\bz}X)$ and $\lbrace\bo\rbrace\times(T_{\bar{\bz}}\overline{X})^{\perp_{\R}}\times(T_{\bz}X)^{\perp_{\R}}$. The image of the derivative of the composition of $\pi_{\V^2}$ with the projection $\pi_{\V}\colon\V^2\to\V$ onto the first component contains both $T_{\bz}X$ and $(T_{\bar{\bz}}\overline{X})^{\perp_{\R}}$. Since $(T_{\bar{\bz}}\overline{X})^{\perp_{\R}}$ is Hermitian orthogonal to $T_{\bz}X$, these spaces intersect trivially and their union span all of $\V$, thus
the derivative of $\pi_{\V}\circ\pi_{\V^2}$ at $(\bz,\bar{\bz},\bz,\bar{\bz})$ is surjective and $\pi_{\V}\circ\pi_{\V^2}$ is dominant.
\end{proof}

\begin{rmk}\label{edcorr}
The \emph{ED correspondence} is defined in a similar way in \cite{dhost2014} by considering the zero locus of the critical ideal \eqref{conded} of the $\ED$ in the polynomial ring $\C[\bz,\bu]$. 
\end{rmk}

The next corollary provide a parametrization of $\vH_X$ when the variety $X$ is parametrized itself.

\begin{Coro}\label{corresp}
If $X$ is (uni-)rational then so is the vHD correspondence $\vH_{X}$.
\end{Coro}
\begin{proof}
Follows as \cite[Corollary 4.2]{dhost2014}. However, we present the proof that provides the formula for the parametrization.

Let $\psi\colon\C^m\to\C^n$ be a rational map that parametrizes $X$, where $m=\dim(X)=n-c$. Its Jacobian $J(\psi)$ is an $n\times m$ matrix of rational functions in the standard coordinates $t_1,\ldots,t_m$ on $\C^m$. The columns of $J(\psi)$ span the tangent space of $X$ at the point $\psi(\bz)$ for generic $\bz\in\C^m$. The left kernel of $J(\psi)$ is a linear space of dimension $c$. We can write down a basis $\lbrace\beta_1(\bz),\ldots,\beta_c(\bz)\rbrace$ of that kernel by applying Cramer’s rule to the matrix $J(\psi)$. In particular, the maps $\beta_k$ for $k=1,\ldots,c$ will also be rational functions in the $z_j$ for $j=1,\ldots,m$. In the end, the map
\begin{align*}
    &\C^{2m}\times\C^{2c}\to\vH_{X}\\
    &(\bz,\bw,\bt,\bs)\mapsto\left(\psi(\bz),\bar{\psi}(\bw),\psi(\bz)+\sum_{k=1}^ct_k\bar{\beta}_k(\bw),\bar{\psi}(\bw)+\sum_{k=1}^cs_k\beta_k(\bz)\right)
\end{align*}
is a parametrization of $\vH_{X}$ that is birational if and only if $\psi$ is birational.
\end{proof}

If $X$ is an affine cone we consider the closure of the image
\begin{align*}
\vH_{X}\cap&((\V\setminus\lbrace\bo\rbrace)^2\times\V^2)\longrightarrow(\Pro\V)^2\times\V^2\\
&((\bz,\bw),(\bu,\bv))\longmapsto(([\bz],[\bw]),(\bu,\bv))
\end{align*}
This closure is called the \emph{projective virtual Hermitian Distance correspondence} (projective vHD correspondence) of $X$, and it is denoted $\vPH_{X}$. It possesses the properties resumed in the following theorem.

\begin{Theo}\label{proj}
Let $X\subseteq\V$ be an irreducible affine cone. The projective vHD correspondence $\vPH_{X}$ is a $2n$ dimensional irreducible variety inside $(\Pro\V)^2\times\V^2$. It is the zero set of the Hermitian critical ideal \eqref{cond3}. The projection $\pi_{X\times\overline{X}}$ is a vector bundle of rank $2c+2$. The projection $\pi_{\V^2}$ has generic fibers of cardinality equal to $\vHD(X)$.
\end{Theo}
\begin{proof}
Follows as \cite[Theorem 4.4]{dhost2014}. The only non trivial point is the fact that $\vPH_{X}$ is the zero set of the Hermitian critical ideal of expression \eqref{cond3}.

Firstly, if $(\bz,\bw,\bu,\bv)\in\vH_{X}$, then $([\bz],[\bw],\bu,\bv)$ lies in the zero locus of the Hermitian critical ideal in \eqref{cond3}. Conversely, if
$([\bz],[\bw],\bu,\bv)$ lies in the variety of that ideal, then there exist unique $\mu_{\bz},\mu_{\bw}\in\C$ such that $(\mu_{\bz}\bz,\mu_{\bw}\bw,\bu,\bv)\in\vH_{X}$. If both $\mu_{\bz}$ and $\mu_{\bw}$ are non zero, then this means that $([\bz],[\bw],\bu,\bv)$ lies in $\vPH_{X}$. Suppose $\mu_{\bz}=0$ and $\mu_{\bw}\neq 0$. In this case
$\bu\perp_{\R}T_{\bw}\overline{X}$ and $\bv-\bw\perp_{\R}T_{\bz}X$, hence $([\varepsilon\bz],[\bw],\bu+\varepsilon\bz,\bv)\in\vPH_{X}$ for any non zero $\varepsilon\in\C$, where we used the fact $T_{\bz}X=T_{\varepsilon\bz}X$. The limit of $([\varepsilon\bz],[\bw],\bu+\varepsilon\bz,\bv)$ for $\varepsilon\rightarrow 0$ equals
$([\bz],[\bw],\bu,\bv)$, so the latter point still lies in the projective ED correspondence $\vPH_{X}$. Similarly arguments apply in general if $\mu_{\bz}\mu_{\bw}=0$.
\end{proof}

Since we are interested in critical points, we define the \emph{Hermitian Distance correspondence} (HD correspondence) $\mathcal{H}_{X}$ as the image of the projection
\begin{align*}
    \vH_{X}\cap\lbrace (\bz,\bar{\bz},\bu,\bar{\bu})\mid(&\bz,\bu)\in\V^{2}\rbrace\to\mathcal{H}_{X}\subseteq\V\times\V\\
    &(\bz,\bw,\bu,\bv)\mapsto(\bz,\bu)
\end{align*}
and similarly the \emph{projective Hermitian Distance correspondence} (projective HD correspondence) $\mathcal{PH}_{X}$ if $X$ is a projective variety.

In particular, if $X$ is parametrized by $\psi$, from Corollary~\ref{corresp}, a parametrization of the HD correspondence is given by the map
\begin{align*}
    \C^{m}&\times\C^{c}\to\mathcal{H}_{X}\subseteq\V\times\V\\
    &(\bz,\bt)\mapsto\left(\psi(\bz),\psi(\bz)+\sum_{k=1}^ct_k\bar{\beta}_k(\bar{\bz})\right)
\end{align*}

We apply these concepts on Example~\ref{parab1}.\newline

We briefly discuss the notion of dual variety. For a more general introduction on dual varieties that does not require an inner product see \cite[Section 1]{gkz1994}.\newline

Consider an affine cone $X\subseteq\V$, or equivalently the corresponding projective variety $X\subseteq\Pro\V$ and let $\V^{\vee}$ be the dual complex vector space of $\V$. Such a variety possesses a dual variety $X^{\vee}\subseteq\V^{\vee}$, which is defined as the following Zariski closure
\begin{equation*}
  X^{\vee}\coloneqq\overline{\lbrace \bs\in\V^{\vee}\mid\exists\ \bz\in X_{\reg}\colon\bs\perp_{\R} T_{\bz}X\rbrace}^{Z}\subseteq\V^{\vee}.
\end{equation*}
Clearly, the conjugate operation and the dual operation commute in the sense that it holds $\overline{X^{\vee}}=\overline{X}^{\vee}$. Moreover note that, even when $X$ is an affine subspace, the Hermitian orthogonal of $X$ is not contained in the dual space $X^{\vee}$. However, it holds the containment $X^{\perp_{\R}}\subseteq X^{\vee}$. We also recall that by the biduality theorem it holds $\left(X^{\vee}\right)^{\vee}=X$, see \cite[Theorem 1.1]{gkz1994}.

A classical known fact is that if we have a bilinear form on $\V$, such as $q(\phantom{\bz},\bar{\phantom{\bz}})$, there exists a canonical isomorphism between $\V$ and $\V^{\vee}$. In particular, we can think of $X^{\vee}$ as a variety in $\V$.\newline

We recall that the \emph{conormal variety} is the Zariski closure
\begin{equation*}
  \mathcal{N}_{X}\coloneqq\overline{\lbrace (\bz,\bs)\in\V^2\mid\bz\in X_{\reg},\ \bs\perp_{\R} T_{\bz}X\rbrace}^{Z}\subseteq\V^2
\end{equation*}
which is the zero set of the bi-homogeneous ideal
\begin{equation*}
    N_{X}\coloneqq \left(I_{X}+\left\langle\text{$(c+1)$-minors of $\begin{bmatrix}
        \bs\\
        J(f)
  \end{bmatrix}$}\right\rangle\right)\colon\left(I_{X_{\sing}}\right)^{\infty}\subseteq\C[\bz,\bs],
\end{equation*}
where $c$ is the codimension of $X$. Recall that we denote $X_{\circ}=X_{\reg}\times\overline{X}_{\reg}$. The following theorem characterizes the action of the dual operator on the Hermitian distance problem.

\begin{Theo}\label{dual}
Let $X\subseteq\V$ be an irreducible affine cone and $(\bu,\bv)\in\V^2$ be a point. The map $(\bz,\bw)\mapsto (\bv-\bw,\bu-\bz)$ is a bijection from the set of solutions of the Hermitian critical ideal of $X$ of $(\bu,\bv)$ to the set of solutions of the Hermitian critical ideal of $X^{\vee}$ of $(\bv,\bu)$ that sends critical points to critical points. In particular, there hold
\begin{equation*}
	\vHD(X)=\vHD(X^{\vee})\qquad\text{and}\qquad\HD(X)=\HD(X^{\vee}).
\end{equation*}
Moreover, the map is proximity reversing for critical points.
\end{Theo}
\begin{proof}
If $(\bz,\bw)\in X_{\circ}$ is in the zero set of the Hermitian critical ideal of $X$ of $(\bu,\bv)$, then $\bv-\bw\perp_{\R}T_{\bz}X$ and hence $\bv-\bw\in X^{\vee}$ or more generally $(\bz,\bv-\bw)\in\mathcal{N}_X$. Similarly, we can say $(\bw,\bu-\bz)\in\mathcal{N}_{\overline{X}}$. Since $(\bu,\bv)$ is generic, we can assume all points $(\bv-\bw,\bu-\bz)$ thus obtained are in $(X^{\vee})_{\circ}$. By the biduality theorem, we have $\bu-(\bu-\bz)=\bz\perp_{\R}T_{\bv-\bw}X^{\vee}$ and similarly $\bv-(\bv-\bw)\perp_{\R} T_{\bu-\bz}\overline{X}^{\vee}$. Thus, $(\bv-\bw,\bu-\bz)\in(X^{\vee})_{\circ}$ is a solution of the Hermitian critical ideal of $X^{\vee}$ of $(\bv,\bu)$.

Considering the case $\bv=\bar{\bu}$ we get the part about critical points.

Applying the same argument to $X^{\vee}$, since $(X^{\vee})^{\vee}=X$, we prove bijectivity.

For the last statement, we restrict our attention to points such that $\bv=\bar{\bu}$ and $\bw=\bar{\bz}$, in particular we compute the distance by considering the first components of the tuples. Observe that it holds $\bu-\bz\perp_{\R}\bw\in T_{\bw}\overline{X}$, or in other words $\langle\bu-\bz,\bw\rangle_{\R}=\langle\bu-\bz,\bar{\bz}\rangle_{\R}=0$, thus $\bu-\bz$ and $\bz$ are Hermitian orthogonal and
\begin{equation*}
    \|\bu-\bz\|_{2}^2+\|\bv-(\bv-\bw)\|_{2}^2=\|\bu-\bz\|_{2}^2+\|\bw\|_{2}^2=\|\bu-\bz\|_{2}^2+\|\bz\|_{2}^2=\|\bu\|_{2}^2.
\end{equation*}
\end{proof}

Introduce now the variables $\bt$ and the ideal
\begin{equation*}
    N_{\overline{X}}\coloneqq \left((I_{X})^{\ast}+\left\langle\text{$(c+1)$-minors of $\begin{bmatrix}
        \bt\\
        J(f)^{\ast}
  \end{bmatrix}$}\right\rangle\right)\colon\left((I_{X_{\sing}})^{\ast}\right)^{\infty}\subseteq\C[\bw,\bt]
\end{equation*}
that is the analogue of the ideal $N_{X}$ for the variety $\overline{X}$. Again following \cite{dhost2014}, duality and the isomorphism above lead us to define the \emph{joint virtual Hermitian Distance correspondence} (joint vHD correspondence) of the cone $X$ in two equivalent ways 
\begin{align*}
  \scalemath{0.9}{\vH_{X}^{\vee}}&\scalemath{0.9}{\coloneqq\overline{\lbrace (\bz,\bw,\bv-\bw,\bu-\bz,\bu,\bv)\in X_{\circ}\times\V^2\times\V^2\mid\bv-\bw\perp_{\R}T_{\bz}X,\ \bu-\bz\perp_{\R}T_{\bw}\overline{X}\rbrace}^{Z}}\\
  &\scalemath{0.9}{=\overline{\lbrace (\bu-\bt,\bv-\bs,\bs,\bt,\bu,\bv)\in\V^2\times (X^{\vee})_{\circ}\times\V^2\mid\bu-\bt\perp_{\R}T_{\bs}X^{\vee},\ \bv-\bs\perp_{\R}T_{\bt}\overline{X}^{\vee}\rbrace}^{Z}},
\end{align*}
and similarly the \emph{projective joint virtual Hermitian Distance correspondence} (projective joint vHD correspondence) $\vPH_{X}^{\vee}$.

Recall from Section~\ref{sec:first} that the variety $X_{\circ}$ is never contained in $\tilde{Q}\coloneqq V\left(\langle\bz,\bw\rangle_{\R}\right)$. With the next result we exploit the connection between the vHD correspondence and the dual of an affine variety. 

\begin{Prop}
Let $X\subseteq\V$ be an irreducible affine cone. The projective joint vHD correspondence $\vPH_{X}^{\vee}$ is an irreducible $2n$-dimensional variety in $\left(\Pro\V\right)^2\times\left(\Pro\V\right)^2\times\V^{2}$. It is the zero set of the hexa-homogeneous ideal
\begin{equation}\label{dualideal}
  \left(N_{X}+N_{\overline{X}}+\left\langle\text{$3$-minors of $\begin{bmatrix}
        \bu\\
        \bt\\
        \bz
    \end{bmatrix}$}\right\rangle+\left\langle\text{$3$-minors of $\begin{bmatrix}
        \bv\\
        \bs\\
        \bw
    \end{bmatrix}$}\right\rangle\right)\colon\langle\langle\bz,\bw\rangle_{\R}\rangle^{\infty}  
\end{equation}
in the polynomial ring $\C[\bz,\bw,\bs,\bt,\bu,\bv]$.
\end{Prop}
\begin{proof}
To see that $\vPH_{X}^{\vee}$ is defined by the ideal \eqref{dualideal}, note first that any point $(\bz,\bw,\bs,\bt,\bu,\bv)$ for which $(\bz,\bw)\in X_{\circ}$ and $\bs\perp_{\R}T_{\bz}X$, $\bt\perp_{\R}T_{\bw}\overline{X}$ satisfies $(\bz,\bs)\in\mathcal{N}_{X}$ and $(\bw,\bt)\in\mathcal{N}_{\overline{X}}$. Now, since $\bz+\bt=\bu$ and $\bw+\bs=\bv$ then the spaces $\spn\lbrace\bz,\bt,\bu\rbrace$ and $\spn\lbrace\bw,\bs,\bv\rbrace$ have dimension at most $2$.

Conversely, let $([\bz],[\bw],[\bs],[\bt],\bu,\bv)$ be in the variety of the ideal \eqref{dualideal}. We can assume $(\bz,\bw)\in X_{\circ}$ and $(\bs,\bt)\in(X^{\vee})_{\circ}$. Since $(\bz,\bs)\in\mathcal{N}_{X}$ and $(\bw,\bt)\in\mathcal{N}_{\overline{X}}$ then $\bs\perp_{\R}T_{\bz}X$ and $\bt\perp_{\R}T_{\bw}\overline{X}$ and in particular $\bs\perp_{\R}\bz$ and $\bt\perp_{\R}\bw$. Applying the hypothesis of saturation must hold that both the pairs $\bs,\bw$ and $\bt,\bz$ are linearly independent, thus there hold $\bu=\mu_1\bz+\lambda_1\bt$ and $\bv=\mu_2\bw+\lambda_2\bs$ for unique constants $\mu_1,\mu_2,\lambda_1,\lambda_2\in\C$. If all those scalars are non zero, then we find that $([\mu_1\bz],[\mu_2\bw],[\lambda_2\bs],[\lambda_1\bt],\bu,\bv)\in\vPH_{X}^{\vee}$. If all are non zero except for $\lambda_2=0$, then the same reasoning as in the proof of Theorem~\ref{proj} proves $([\mu_1\bz],[\mu_2\bw],[\varepsilon\bs],[\lambda_1\bt],\bu,\bv+\varepsilon\bs)\in\vPH_{X}^{\vee}$ for any non zero $\varepsilon\in\C$ and the limit point $([\mu_1\bz],[\mu_2\bw],[\bs],[\lambda_1\bt],\bu,\bv)$ for $\varepsilon\rightarrow 0$ lies in $\vPH_{X}^{\vee}$.

Similar arguments apply in general when $\mu_1\mu_2\lambda_1\lambda_2=0$.
\end{proof}

In our context could be more natural to consider a different isomorphism that we now mention. In the presence of a sesquilinear form, that is $q$, there exists a canonical isomorphism between $\V$ and $\overline{\V^{\vee}}=\overline{\V}^{\vee}$ and in particular it holds
\begin{equation*}
  \overline{X}^{\vee}=\overline{\lbrace \bs\in\V\mid\exists\ \bz\in X_{\reg}\colon\bar{\bs}\perp_{\R} T_{\bz}X\rbrace}^{Z}\subseteq\V.
\end{equation*}


\section{The HD discriminant}\label{sec:hddisc}\thispagestyle{plain}

For the Euclidean distance problem, the \emph{ED discriminant} $\Sigma_X\subseteq\V$ of a real algebraic variety $X\subseteq\V$ is defined in \cite{dhost2014}. This variety possesses many properties, in particular the number of critical points is constant on the different chambers defined by $\Sigma_X$.

We apply similar concepts for the Hermitian case.\newline

Let $X\subseteq\V$ be an algebraic variety, by Theorem~\ref{corr} the $\vHD(X)$ is the cardinality of the generic fiber of the projection $\pi_{\V^2}\colon\vH_{X}\subseteq X\times\overline{X}\times\V^2\rightarrow\V^{2}$ of the vHD correspondence, see Subsection~\ref{ssec:corrdual}.

\begin{Def}
Let $X\subseteq\V$ be an algebraic variety, the branch locus of $\pi_{\V^2}$ defines a variety $\vX_{X}\subseteq\V^2$ that we call the \emph{virtual Hermitian Distance discriminant} (vHD discriminant) of $X$.
\end{Def}

The Nagata-Zariski purity theorem implies that the last definition is well posed and that the vHD discriminant is typically a hypersurface. We are interested in its defining polynomial which can be interpreted as a polynomial in $\C[\bu,\bv]$. In particular, this polynomial is invariant under the action of the restricted map $\ast$ defined in Section~\ref{sec:first}.

Moreover, the vHD discriminant is contained in the Zariski closure 
\begin{equation*}
    \overline{\lbrace (\bu,\bv)\in\V^2\mid\#\pi_{\V^2}^{-1}(\bu,\bv)<\vHD(X)\rbrace}^{Z}
\end{equation*}
where $\#$ denotes the set theoretical cardinality. In particular, $\vX_{X}$ contains the points $(\bu,\bv)$ in which a solution of the Hermitian critical ideal of $(\bu,\bv)$ has multiplicity greater than one.

If $X\subseteq\Pro\V$ is projective, since the variety $\vH_{X}$ is defined by tetra-homogeneous equations in $\bz,\bw,\bu,\bv$, also the branch locus is defined by homogeneous equations.

Since we are interested in critical points, we specialize our definition.

\begin{Def}
Let $X\subseteq\V$ be an algebraic variety, the \emph{Hermitian Distance discriminant} (HD discriminant) of $X$ is the set $\Xi_{X}\coloneqq\lbrace\bu\in\V\mid(\bu,\bar{\bu})\in\vX_{X}\rbrace\subseteq\V$.
\end{Def} 

Note that the HD discriminant $\Xi_{X}$ is defined via the generalized polynomial obtained by setting $\bv=\bar{\bu}$ in the equation defining the vHD discriminant $\vX_{X}$ and thus, in general, is not a complex algebraic variety. However, it is a real algebraic variety. 

On the other hand, the real locus of $\Xi_{X}$ coincides with the projection onto the first component of the real locus of $\vX_{X}$. In other words, the defining polynomial of the real locus of $\Xi_{X}$ is obtained setting $\bv=\bu$ in the polynomial defining $\vX_X$ and thus is a real algebraic variety contained in $\Xi_{X}$.

The following result directly follows from the last part of Proposition~\ref{HDvsED}.

\begin{Prop}
Let $X\subseteq\V$ be a real algebraic variety. The real locus of $\Sigma_X\subseteq\V$ is contained in the real locus of $\Xi_{X}\subseteq\V$.
\end{Prop}
\begin{proof}
Follows setting $\bv=\bu$ real and $\bw=\bz$, and noting that with these choices the two ideals defining the HD correspondence and the ED correspondence coincide, see Remark~\ref{edcorr}.
\end{proof}

For the sake of simplicity, we will use the symbol $\vX_{X}$ to also intend the equation defining the vHD discriminant and similarly for $\Xi_{X}$.\newline

We can consider three cases for the points in the vHD discriminant:
\begin{enumerate}[\normalfont i)]
  \item There is exactly one solution of the critical ideal with multiplicity equal to $2$.
  \item There are two or more solutions of the critical ideal with multiplicity equal to $2$.
  \item Any other possibility.
\end{enumerate}
In case i), setting $\bv=\bar{\bu}$, the number of critical points in a neighborhood could change by $2$. On the other hand, in cases ii) and iii), setting $\bv=\bar{\bu}$, the jump between the numbers of critical points could be higher. In particular, the multiplicity of at least a solution of the critical ideal is greater than $2$ for the points of case iii), moreover these points are contained in the singular locus of $\vX_X$. By a slight abuse of notation, we call the singular locus of the HD discriminant $\Xi_{X}$ the set obtained by setting $\bv=\bar{\bu}$ in the singular locus of $\vX_{X}$.

We have already seen in Example~\ref{exnotcon} a variety for which there are jumps by more than $2$ between the numbers of $\HD$. More generally, let $X\subseteq\V$ and $Y\subseteq\W$ be algebraic varieties and consider $X\times Y$. The critical ideal of $X\times Y$ is the sum of the critical ideals of $X$ and $Y$. The HD discriminant $\Xi_{X\times Y}$ is $(\Xi_{X}\times\W)\cup(\V\times\Xi_{Y})\subseteq\V\times\W$ and if $\min\lbrace\vHD(X),\vHD(Y)\rbrace>1$ any of its points satisfy condition ii) or iii) above.\newline

From the results on the ED discriminant $\Sigma_X$ and the discussion above, we obtain the following.

\begin{Prop}
Let $X\subseteq\V$ be a real algebraic variety. The set $\Xi_{X}\setminus\Sigma_{X}\subseteq\V$ contains a dense subset of the set of points that admit only non real critical points of the Hermitian distance with multiplicity greater than one. The real locus of the variety $\Sigma_{X}\subseteq\V$ contains a dense subset of the set of points that admit only real critical point of the Hermitian distance with multiplicity greater than one.\newline
\end{Prop}

The next result follows from the properties of the Hermitian Killing form, in particular see \cite[Corollary 3.14]{f2024}.

\begin{Coro}
Let $X$ be an algebraic variety, the set $\Xi_X$ is contained in the zero locus of the determinant of the Hermitian Killing form of the system given by the Hermitian critical ideal \eqref{cond}. Moreover, the number of critical points changes on a path crossing $\Xi_X$ if the sign of the evaluation of $\Xi_X$ changes.
\end{Coro}
\begin{proof}
The zero locus of the determinant of the Hermitian Killing form of the system given by the Hermitian critical ideal determines exactly the points $(\bu,\bv)$ for which a solution of the Hermitian critical ideal has multiplicity greater than one. The latter is the definition of the vHD discriminant $\vX_X$ restricted to points such that $\bv=\bar{\bu}$.
\end{proof}


\subsection{Complex evolute}\label{ssec:evolute}

Here we present known properties of the ED discriminant of a real algebraic curve. We then apply similar ideas to the case of the HD discriminant of an algebraic curve.\newline

\begin{figure}[H]
\centering
\includegraphics[height=5cm, width=5cm]{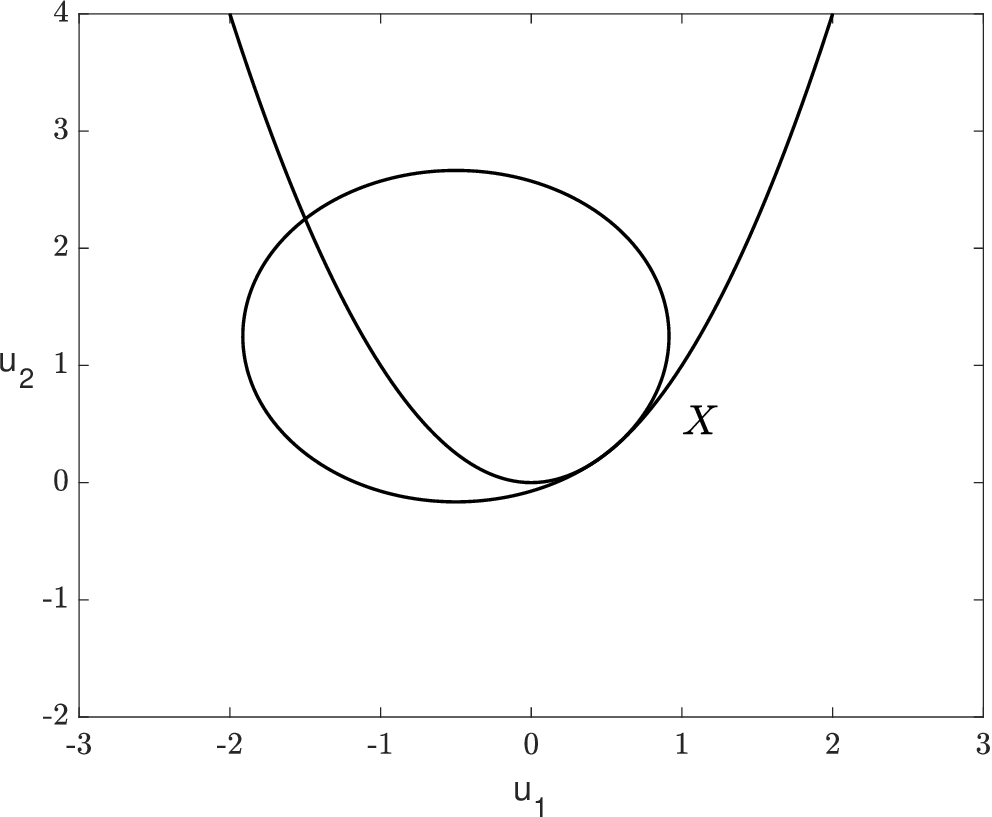}
\qquad
\includegraphics[height=5cm, width=5cm]{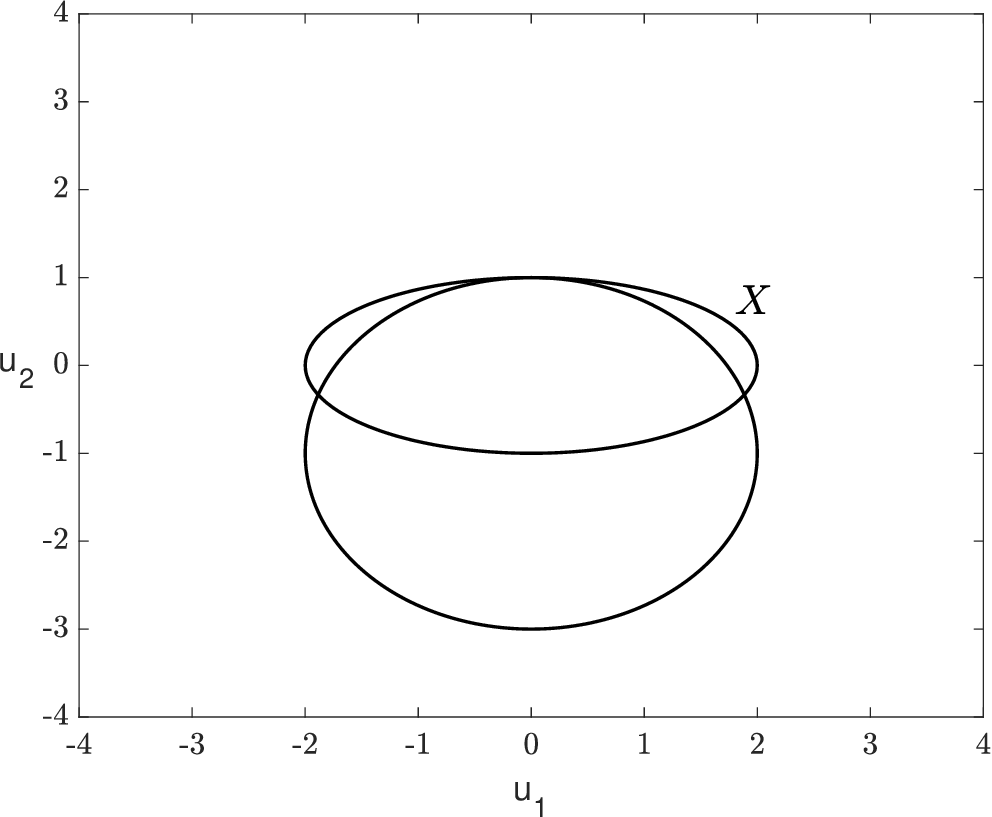}
\caption{Osculating circles of the parabola $X=V(z_2-z_1^2)$ at the point $\bu=(1/2,1/4)$ (Left) and of the ellipse $X=V(z_1^2+4z_2^2-4)$ at the point $\bu=(0,1)$ (Right).}\label{oscfig}
\end{figure}

Let us recall that an \emph{osculating circle} of a curve $X\subseteq\R^2$ is a circle that is tangent to the curve at a regular point and has the same curvature in the same direction as the curve at the point of tangency. Two simple examples of osculating circles are shown in Figure~\ref{oscfig}.

We are ready to state the definition of the evolute of a real algebraic curve, a rigorous introduction can be found in \cite{bg1992}, other works about the more general concept of \emph{Focal loci} are \cite{c1998,ct2007}.

\begin{Def}\label{defevo}
Let $X=V(f)\subseteq\R^2$ with $f\in\R[\bz]$ be a plane curve, the \emph{evolute} of $X$ is the union of the centers of the osculating circles of $X$. Its defining polynomial is the generator of the ideal obtained by eliminating the variables $\bz$ from the ideal
\begin{equation*}
  I_{X}^{E}\coloneqq\left\langle f, g, \varphi \right\rangle\colon\langle\partial_{z_1}f,\partial_{z_2}f\rangle^{\infty}\subseteq\C[\bz,\bu],
\end{equation*}
where $g\coloneqq(u_2-z_2)\partial_{z_1}f-(u_1-z_1)\partial_{z_2}f$ and $\varphi\coloneqq(\partial_{z_1}g)(\partial_{z_2}f)-(\partial_{z_2}g)(\partial_{z_1}f)$.
\end{Def}

In particular, the evolute coincides with the envelope of the normal lines of $X$, see \cite{bg1992}.

Moreover, if $X\subseteq\R^2$ is a smooth curve parametrized by $\psi=(\psi_1,\psi_2)$, the evolute can be obtained using the formula 
\begin{equation}\label{evolute}
    \psi^{E}\coloneqq\left(\psi_1-\psi_2^{\prime}\frac{(\psi_1^{\prime})^2+(\psi_2^{\prime})^2}{\psi_1^{\prime}\psi_2^{\prime\prime}-\psi_1^{\prime\prime}\psi_2^{\prime}},\ \psi_2+\psi_1^{\prime}\frac{(\psi_1^{\prime})^2+(\psi_2^{\prime})^2}{\psi_1^{\prime}\psi_2^{\prime\prime}-\psi_1^{\prime\prime}\psi_2^{\prime}}\right)
\end{equation}
which plots the center of the osculating circle at $\psi$. It is a matter of computations to check $\langle\psi^{\prime},(\psi^{E})^{\prime}\rangle_{\R}=0$.\newline

The following known result completely solves the problem of finding the ED discriminant of a real curve. The equivalence simply follows by noting that the points of the evolute are precisely those points with a critical point $\bz\in X$ of multiplicity greater than one, or in other terms those points which admit critical points of $g$ on $X$. However, we give a rigorous proof that we will generalize in our setting. 

\begin{Prop}
Let $X=V(f)$ with $f\in\R[\bz]$ be a smooth plane curve. The evolute and the ED discriminant coincide.
\end{Prop}
\begin{proof}
Let $g$ and $\varphi$ be the polynomials of Definition~\ref{defevo}. The critical ideal \eqref{conded} of the $\ED$ is $\langle f,g\rangle\subseteq\C[\bz,\bu]$. Since the definitions of the curves involved in the assertion are local, by the implicit function theorem we can assume $z_2$ is a function of $z_1$ and eliminate $f$. Thus, applying the implicit function theorem, the ED discriminant can be defined as the common zero locus of the polynomials $g$ and
\begin{equation*}
   \partial_{z_1}g=\begin{bmatrix}
       1 & \partial_{z_1}z_2
   \end{bmatrix}\begin{bmatrix}
       \partial_{z_1}g\\
       \partial_{z_2}g
   \end{bmatrix}=\partial_{z_1}g+(\partial_{z_1}z_2)(\partial_{z_2}g)=\partial_{z_1}g-\frac{\partial_{z_1}f}{\partial_{z_2}f}\partial_{z_2}g=\frac{\varphi}{\partial_{z_2}f},
\end{equation*}
Resuming, the ED discriminant can be obtained through the same ideal $I_{X}^{E}$ that defines the evolute.
\end{proof}

The following object is the Hermitian analogue of Definition~\ref{defevo} after the introduction of the variables $\bw$ and $\bv$.

\begin{Def}\label{compevo}
Let $X=V(f)$ be an algebraic curve, we call the \emph{virtual complex evolute} of $X$ the zero locus of the generator of the ideal obtained by eliminating the variables $\bz$ and $\bw$ from the ideal
\begin{equation*}
  I_{X}^{\tilde{E}}\coloneqq\left\langle f,f^{\ast},\tilde{g},\tilde{g}^{\ast},\xi\right\rangle\colon\langle\partial_{z_1}f,\partial_{z_2}f,\partial_{w_1}f^{\ast},\partial_{w_2}f^{\ast}\rangle^{\infty}\subseteq\C[\bz,\bw,\bu,\bv],
\end{equation*}
where $\tilde{g}\coloneqq(v_2-w_2)\partial_{z_1}f-(v_1-w_1)\partial_{z_2}f$ and 
\begin{equation}\label{xi}
    \xi\coloneqq\tilde{\varphi}\tilde{\varphi}^{\ast}-\langle\nabla_{\bz}f,\nabla_{\bw}f^{\ast}\rangle_{\R}^2
\end{equation}
with $\tilde{\varphi}\coloneqq(\partial_{z_2}f)(\partial_{z_1}\tilde{g})-(\partial_{z_1}f)(\partial_{z_2}\tilde{g})$. We call the \emph{complex evolute} of $X$ the zero locus of the generalized polynomial obtained setting $\bv=\bar{\bu}$.
\end{Def}

Note that, when considering points such that $\bv=\bar{\bu}$ and critical points $\bw=\bar{\bz}$ the formula simplifies to
\begin{equation*}
    \xi=|\tilde{\varphi}|^2-\|\nabla_{\bz}f\|_{2}^4=\left(|\tilde{\varphi}|-\|\nabla_{\bz}f\|_{2}^2\right)\left(|\tilde{\varphi}|+\|\nabla_{\bz}f\|_{2}^2\right)
\end{equation*}
and in particular we can bound from above the two terms as
\begin{equation*}
    |\tilde{\varphi}|\pm\|\nabla_{\bz}f\|_{2}^2\leq\|\nabla_{\bz}\tilde{g}\|_{2}\|\nabla{\bz}f\|_{2}\pm\|\nabla_{\bz}f\|_{2}^2=\|\nabla_{\bz}f\|_{2}\left(\|\nabla_{\bz}\tilde{g}\|_{2}\pm\|\nabla_{\bz}f\|_{2}\right).\newline
\end{equation*}

We now present the connection between the HD discriminant and the complex evolute of a curve.

\begin{Prop}
Let $X=V(f)$ be a smooth algebraic curve. The virtual complex evolute of $X$ and the vHD discriminant $\vX_{X}$ coincide. In particular, the complex evolute of $X$ and the HD discriminant $\Xi_{X}$ coincide and this is exactly the set where the number of critical points could change in a neighborhood.
\end{Prop}
\begin{proof}
Let $\tilde{g}$ be the polynomial of Definition~\ref{compevo}. The Hermitian critical ideal takes the form $\left\langle f,f^{\ast},\tilde{g},\tilde{g}^{\ast}\right\rangle$ and we can assume $z_2$ and $w_2$ to be holomorphic functions of $z_1$ and $w_1$ respectively and thus eliminate $f$ and $f^{\ast}$. Moreover, we can assume $w_1$ to be a holomorphic function of $z_1$ and eliminate $\tilde{g}^{\ast}$. Thus, applying the implicit function theorem to derive $\partial_{z_1}z_2(z_1)$ and $\partial_{w_1}w_2(w_1)$ and the chain rule, the HD discriminant is the common zero locus of the polynomials $\tilde{g}$ and
\begin{align*}
    \partial_{z_1}\tilde{g}&=\begin{bmatrix}
       1 & \partial_{z_1}z_2 & \partial_{z_1}w_1 & \partial_{z_1}w_2
   \end{bmatrix}\begin{bmatrix}
       \partial_{z_1}\tilde{g}\\
       \partial_{z_2}\tilde{g}\\
       \partial_{w_1}\tilde{g}\\
       \partial_{w_2}\tilde{g}
   \end{bmatrix}\\
   &=\begin{bmatrix}
       1 & -\frac{\partial_{z_1}f}{\partial_{z_2}f} & \partial_{z_1}w_1 & -\partial_{z_1}w_1\frac{\partial_{w_1}f^{\ast}}{\partial_{w_2}f^{\ast}}
   \end{bmatrix}\begin{bmatrix}
       \partial_{z_1}\tilde{g}\\
       \partial_{z_2}\tilde{g}\\
       \partial_{w_1}\tilde{g}\\
       \partial_{w_2}\tilde{g}
   \end{bmatrix}\\
   &=\partial_{z_1}\tilde{g}-\frac{\partial_{z_1}f}{\partial_{z_2}f}\partial_{z_2}\tilde{g}+\partial_{z_1}w_1\left(\partial_{w_1}\tilde{g}-\frac{\partial_{w_1}f^{\ast}}{\partial_{w_2}f^{\ast}}\partial_{w_2}\tilde{g}\right)\\
   &=\frac{(\partial_{z_2}f)(\partial_{z_1}\tilde{g})-(\partial_{z_1}f)(\partial_{z_2}\tilde{g})+\frac{\partial_{z_2}f}{\partial_{w_2}f^{\ast}}\partial_{z_1}w_1\left((\partial_{z_2}f)(\partial_{w_2}f^{\ast})+(\partial_{z_1}f)(\partial_{w_1}f^{\ast})\right)}{\partial_{z_2}f}.
\end{align*}
Now, we use again the implicit function theorem to explicit $\partial_{z_1}w_1(z_1)$. Starting from the equalities
\begin{align*}
    \partial_{z_1}w_1=-\frac{\partial_{z_1}\tilde{g}^{\ast}}{\partial_{w_1}\tilde{g}^{\ast}}&=-\frac{\partial_{z_1}\left((u_2-z_2)\partial_{w_1}f^{\ast}-(u_1-z_1)\partial_{w_2}f^{\ast}\right)}{\partial_{w_1}\left((u_2-z_2)\partial_{w_1}f^{\ast}-(u_1-z_1)\partial_{w_2}f^{\ast}\right)}\\
    &=-\frac{1}{\partial_{z_2}f}\frac{(\partial_{z_1}f)(\partial_{w_1}f^{\ast})+(\partial_{z_2}f)(\partial_{w_2}f^{\ast})}{(u_2-z_2)\partial_{w_1}\left(\partial_{w_1}f^{\ast}\right)-(u_1-z_1)\partial_{w_1}\left(\partial_{w_2}f^{\ast}\right)},
\end{align*}
we apply the equalities
\begin{align*}
    &\partial_{w_1}(\partial_{w_1}f^{\ast})=\begin{bmatrix}
        1 & -\frac{\partial_{w_1}f^{\ast}}{\partial_{w_2}f^{\ast}}
    \end{bmatrix}\begin{bmatrix}
        \partial_{w_1}^2f^{\ast}\\
        \partial_{w_1}\partial_{w_2}f^{\ast}
    \end{bmatrix}=\partial_{w_1}^2f^{\ast}-\frac{\partial_{w_1}f^{\ast}}{\partial_{w_2}f^{\ast}}\partial_{w_1}\partial_{w_2}f^{\ast}\\
    &\partial_{w_1}(\partial_{w_2}f^{\ast})=\begin{bmatrix}
        1 & -\frac{\partial_{w_1}f^{\ast}}{\partial_{w_2}f^{\ast}}
    \end{bmatrix}\begin{bmatrix}
        \partial_{w_1}\partial_{w_2}f^{\ast}\\
        \partial_{w_2}^2f^{\ast}
    \end{bmatrix}=\partial_{w_1}\partial_{w_2}f^{\ast}-\frac{\partial_{w_1}f^{\ast}}{\partial_{w_2}f^{\ast}}\partial_{w_2}^2f^{\ast}
\end{align*}
obtained by the chain rule to write
\begin{equation*}
    \partial_{z_1}w_1=-\frac{\partial_{w_2}f^{\ast}}{\partial_{z_2}f}\frac{(\partial_{z_2}f)(\partial_{w_2}f^{\ast})+(\partial_{z_1}f)(\partial_{w_1}f^{\ast})}{(\partial_{w_2}f^{\ast})(\partial_{w_1}\tilde{g}^{\ast})-(\partial_{w_1}f^{\ast})(\partial_{w_2}\tilde{g}^{\ast})}.
\end{equation*}
Now, applying this substitution in the formula for $\partial_{z_1}\tilde{g}$ results $\xi$ of equation \eqref{xi} in the numerator, thus we get the first part of the statement.

The second statement follows from the definitions of the complex evolute and critical points.
\end{proof}

In the following subsections we consider $X$ to be a real algebraic curve and we will provide a useful description of the HD discriminant. 


\subsection{Outward evolute}

\begin{figure}[H]
\centering
\includegraphics[height=5cm, width=5cm]{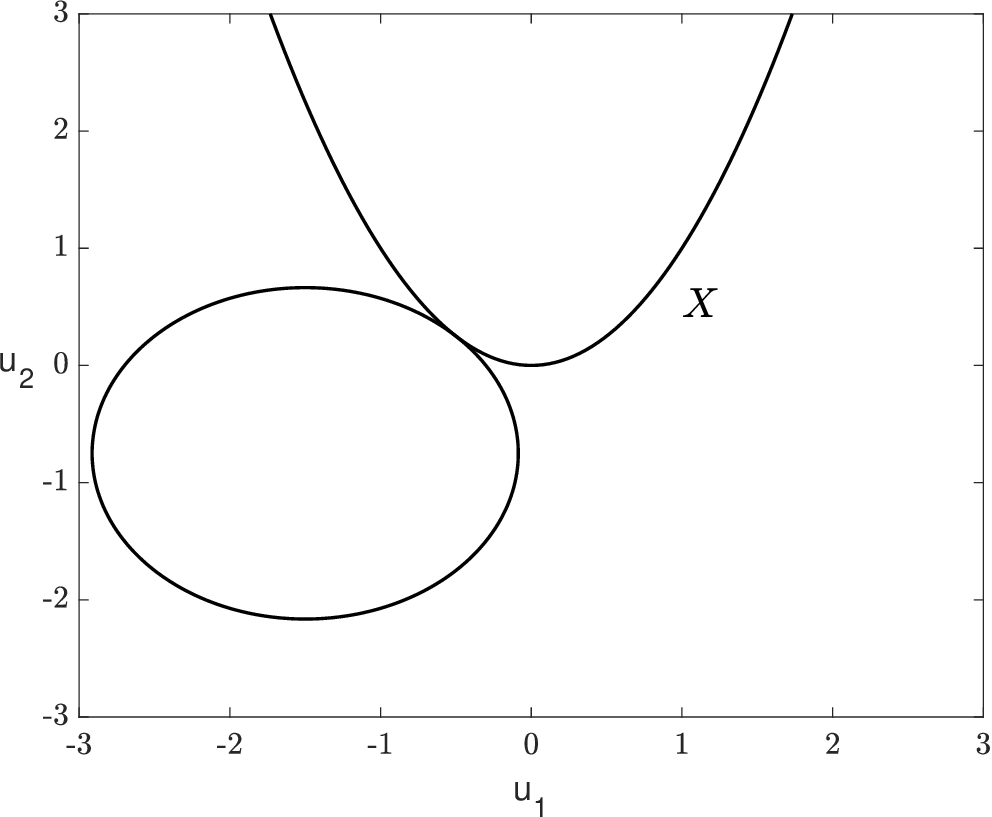}
\qquad
\includegraphics[height=5cm, width=5cm]{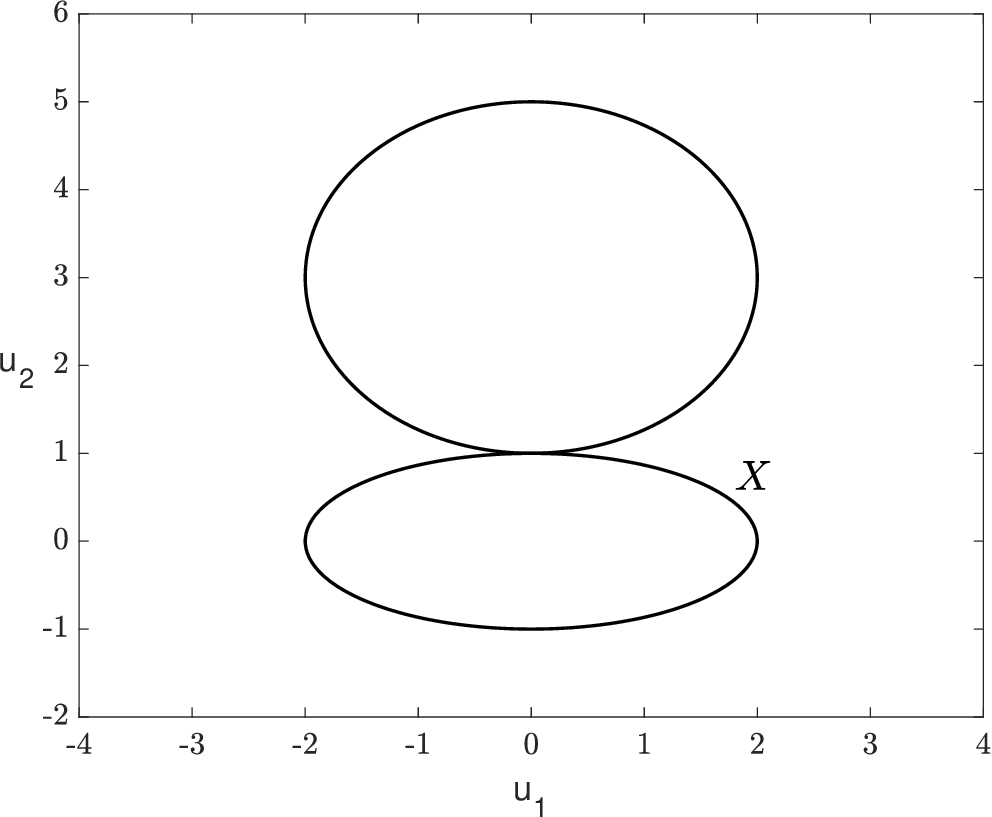}
\caption{Outward osculating circles of the parabola $X=V(z_2-z_1^2)$ at the point $\bu=(-1/2,1/4)$ (Left) and of the ellipse $X=V(z_1^2+4z_2^2-4)$ at the point $\bu=(0,1)$ (Right).}\label{outoscfig}
\end{figure}

To our purpose, in this subsection it is useful to introduce the notion of an \emph{outward osculating circle} of a curve $X\subseteq\R^2$, that is a circle that is tangent to the curve at a regular point and has the same curvature with opposite direction as the curve at point of tangency. Two simple examples of outward osculating circles are shown in Figure~\ref{outoscfig}.

Using the outward osculating circles we define our object of interest. In the following definition and in the next result, for a vector $\bz$, we will simply write $\|\bz\|_{\R}^2$ to denote $\langle\bz,\bz\rangle_{\R}$.

\begin{Def}
Let $X=V(f)\subseteq\R^2$ with $f\in\R[\bz]$ be a plane curve, the \emph{outward evolute} of $X$ is the union of the centers of the outward osculating circles of $X$. Its defining polynomial is the generator of the ideal obtained eliminating the variables $\bz$ from the ideal
\begin{equation*}
  I_{X}^{\EE}\coloneqq\left\langle f, g, \varphi-2\|\nabla_{\bz}f\|_{\R}^2\right\rangle\colon\langle\partial_{z_1}f,\partial_{z_2}f\rangle^{\infty}\subseteq\C[\bz,\bu],
\end{equation*}
where $g$ and $\varphi$ are the polynomials of Definition~\ref{defevo}.
\end{Def}

If $\psi$ is a parametrization of the smooth curve $X\subseteq\R^2$ we can modify formula \eqref{evolute} of the evolute as 
\begin{equation*}
    \psi^{\EE}\coloneqq\left(\psi_1+\psi_2^{\prime}\frac{(\psi_1^{\prime})^2+(\psi_2^{\prime})^2}{\psi_1^{\prime}\psi_2^{\prime\prime}-\psi_1^{\prime\prime}\psi_2^{\prime}},\ \psi_2-\psi_1^{\prime}\frac{(\psi_1^{\prime})^2+(\psi_2^{\prime})^2}{\psi_1^{\prime}\psi_2^{\prime\prime}-\psi_1^{\prime\prime}\psi_2^{\prime}}\right)
\end{equation*}
to obtain the outward evolute of $X$. In particular, $\psi^{\EE}=2\psi-\psi^{E}$ and $\langle\psi^{\prime},(\psi^{\EE})^{\prime}\rangle_{\R}=2\|\psi^{\prime}\|_{\R}^2$. Moreover, $\psi-\psi^{E}=\psi^{\EE}-\psi$ is the radial curve of $X$.\newline

The following proposition shed some light on the real locus of the complex evolute of a real algebraic curve defined in the previous subsection.

\begin{Prop}\label{compev}
Let $X=V(f)$ with $f\in\R[\bz]$ be a smooth algebraic curve. The union of the evolute and the outward evolute is the real locus of the complex evolute. Moreover, for points of
the evolute at least two real critical points coincide, while for points of the outward evolute at least two non real critical points coincide.
\end{Prop}
\begin{proof}
Let be $\bv=\bu$ and write
\begin{equation*}
        \tilde{\varphi}=(\partial_{z_2}f)(\partial_{z_1}\tilde{g})-(\partial_{z_1}f)(\partial_{z_2}\tilde{g})=\varphi+(\partial_{z_2}f)(\partial_{z_1}(\tilde{g}-g))-(\partial_{z_1}f)(\partial_{z_2}(\tilde{g}-g)).
\end{equation*}
We compute
\begin{align*}
    \partial_{z_1}(\tilde{g}-g)&=\partial_{z_1}\left((z_2-w_2)\partial_{z_1}f-(z_1-w_1)\partial_{z_2}f\right)\\
    &=(z_2-w_2)\partial_{z_1}^2f-\partial_{z_2}f-(z_1-w_1)\partial_{z_1}\partial_{z_2}f
\end{align*}
and similarly
\begin{equation*}
    \partial_{z_2}(\tilde{g}-g)=\partial_{z_1}f+(z_2-w_2)\partial_{z_1}\partial_{z_2}f-(z_1-w_1)\partial_{z_2}^2f.
\end{equation*}
For points of the complex evolute at least two solutions of the Hermitian critical ideal of $(\bu,\bu)$ coincide. Moreover, if $\bu$ is real then at least two solutions must coincide at a real critical point, we discussed in Remark~\ref{realpoints} the nature of the solutions in this case. In particular, we can set $\bw=\bz$. In the end, substituting the two equations above, we can write $\tilde{\varphi}=\varphi-\|\nabla_{\bz}f\|_{\R}^2$ and the proposition follows from
\begin{equation*}
    \xi=\tilde{\varphi}\tilde{\varphi}^{\ast}-\langle\nabla_{\bz}f,\nabla_{\bw}f^{\ast}\rangle_{\R}^2=\left(\tilde{\varphi}+\|\nabla_{\bz}f\|_{\R}^2\right)\left(\tilde{\varphi}-\|\nabla_{\bz}f\|_{\R}^2\right)=\varphi\left(\varphi-2\|\nabla_{\bz}f\|_{\R}^2\right).
\end{equation*}
\end{proof}

\begin{Coro}\label{singdisc}
Let $X=V(f)$ with $f\in\R[\bz]$ be an algebraic curve. The outward evolute is contained in the real singular locus of the HD discriminant.
\end{Coro}
\begin{proof}
As discussed in Remark~\ref{realpoints}, when $\bv=\bu$ is real, by conjugation, any non real solution of the critical ideal yields at least one other different solution, and the two solutions form an associated pair. On the outward evolute, at least two points of an associated pair coincide and the coincident point must be a real critical point. In particular, three solutions of the critical ideal coincide thus getting a solution with multiplicity at least $3$ of the critical ideal.
\end{proof}

We now present some examples for the three different cases of non degenerate conic.

\begin{ex}\label{parab1}
(Parabola). The parabola 
\begin{equation*}
    X=V(z_2-z_1^2)\subseteq\C^2
\end{equation*}
satisfies $\vHD(X)=5$ as predicted by Proposition~\ref{degparam} and $\HD(X)=\lbrace 1,3\rbrace$. The value of $\vHD(X)$ can be seen by computing the degree of the critical ideal of $X$ for general parameters $(\bu,\bv)$. This ideal is generated by the polynomials of the system
\begin{equation*}
    \begin{cases}
        p_1(\bz,\bw)=z_2-z_1^2=0\\
        \bar{p}_1(\bw,\bz)=w_2-w_1^2=0\\
        p_2(\bz,\bw)=z_1-u_1+2w_1(z_2-u_2)=0\\ 
        \bar{p}_2(\bw,\bz)=w_1-v_1+2z_1(w_2-v_2)=0  
    \end{cases}
\end{equation*}
where the upper two polynomials are given by the definition of $X$ and the bottom two polynomials are given by the determinant of the $2\times 2$ matrix obtained from the perpendicularity condition. In particular, it follows the containment $\HD(X)\subseteq\lbrace 1,3,5\rbrace$. We now compute the true values of $\HD(X)$.

Set $\bv=\bar{\bu}$. Consider the one-to-one parametrization of the parabola 
\begin{align*}
    \psi\colon&\C\to X\subseteq\C^2\\
    &z\mapsto (z,z^2)
\end{align*}
thus critical points of the Hermitian distance satisfy the equation
\begin{align*}
    \partial_{z}\|\psi(z)-\bu\|_{2}^2&=\partial_z\left[(z-u_1)(\bar{z}-\bar{u}_1)+(z^2-u_2)(\bar{z}^2-\bar{u}_2)\right]\\
    &=(\bar{z}-\bar{u}_1)+2z(\bar{z}^2-\bar{u}_2)=0.
\end{align*}
The system we obtain introducing $\bw$ is
\begin{equation*}
    \begin{cases}
        p(z,w)=(z-u_1)+2w(z^2-u_2)=0\\
        \bar{p}(w,z)=(w-\bar{u}_1)+2z(w^2-\bar{u}_2)=0
    \end{cases}
\end{equation*} 
and it is equivalent to applying substitution to the Hermitian critical ideal and setting $\bv=\bar{\bu}$. The matrix representing the Hermitian Killing form, see \cite{f2024}, with respect to the basis 
\begin{equation*}
    \lbrace [1],[z],[w],[zw],[w^2]\rbrace
\end{equation*}
of the quotient $\C[z,w]/\langle p,p^{\ast}\rangle$ is
\begin{equation*}
    \mathcal{K}_{\C}^1=\begin{bmatrix}
        5 & -\frac{u_1}{2u_2} & -\frac{\bar{u}_1}{2\bar{u}_2} & -2 & \frac{16|u_2|^2\bar{u}_2+\bar{u}_1^2}{4\bar{u}_2^2}\\
        \ast & -2 & \frac{16|u_2|^2\bar{u}_2+\bar{u}_1^2}{4\bar{u}_2^2} & \frac{8u_1\bar{u}_2+\bar{u}_1}{4\bar{u}_2} & \frac{-\bar{u}_1^3-6u_1\bar{u}_2^2}{8\bar{u}_2^2}\\
        \ast & \ast & -2 & \frac{8\bar{u}_1u_2+u_1}{4u_2} & \frac{8u_1\bar{u}_2+\bar{u}_1}{4\bar{u}_2}\\
        \ast & \ast & \ast & 4|u_2|^2+1 & \frac{-2|u_1|^2\bar{u}_2-32|u_2|^2\bar{u}_2-\bar{u}_1^2}{8\bar{u}_2^2}\\
        \ast & \ast & \ast & \ast & 4|u_2|^2+1
    \end{bmatrix}
\end{equation*}
where we assumed without loss of generality $u_2\neq 0$. In particular, the third diagonal entry is negative, thus this matrix possesses at least one negative eigenvalue and from \cite[Corollary 3.14]{f2024} we obtain $\max\HD(X)\leq 4-1=3$. By testing points $\bu$ in the real plane, we get both one or three critical points and thus it follows $\HD(X)=\lbrace 1,3\rbrace$.

The parametrization for the vHD correspondence, see Subsection~\ref{ssec:corrdual}, is
\begin{align*}
  &\C^2\times\C^2\to\mathcal{H}_{X}\subseteq\C^4\times\C^4\\
  &(z,w,s,t)\mapsto\left((z,z^2,w,w^2),(z+2sw,z^2-s,w+2tz,w^2-t)\right)
\end{align*}
Thus, the HD correspondence can be parametrized by
\begin{align*}
  \C&\times\C\to\mathcal{H}_{X}\subseteq\C^2\times\C^2\\
  (&z,s)\mapsto((z,z^2),(z+2s\bar{z},z^2-s))
\end{align*}

The complex evolute of $X$ is the zero locus of
\begin{align*}
  &\scalemath{0.9}{65536|u_2|^{12}+6|u_1|^4(u_1^2\bar{u}_2+\bar{u}_1^2u_2)-|u_2|^4(49152|u_2|^4+10752|u_1|^2)(u_1^2\bar{u}_2+\bar{u}_1^2u_2+|u_2|^2)}\\
  &\scalemath{0.9}{-2|u_1|^6+6|u_1|^4|u_2|^2+384|u_2|^2(|u_1|^4\bar{u}_1^2u_2+|u_1|^4u_1^2\bar{u}_2+u_1^4\bar{u}_2^2+\bar{u}_1^4u_2^2+|u_2|^2u_1^2\bar{u}_2+|u_2|^2\bar{u}_1^2u_2)}\\
  &\scalemath{0.9}{-5952|u_1u_2|^4-1024(u_1^6\bar{u}_2^3+\bar{u}_1^6u_2^3+|u_2|^6)+|u_2|^2(96|u_1|^2-44544|u_2|^4)(\bar{u}_1^2u_2+u_1^2\bar{u}_2+|u_1|^4)}\\
  &\scalemath{0.9}{-172032|u_1|^2|u_2|^8-27(u_1^4\bar{u}_2^2+\bar{u}_1^4u_2^2+|u_1|^8)+(12288|u_2|^4-384|u_1|^2)(u_1^4\bar{u}_2^2+\bar{u}_1^4u_2^2+|u_2|^4)}.
\end{align*}
If we restrict our attention to real points it simplifies to
\begin{equation*}
    (2(2u_2-1)^3-27u_1^2)(8u_2^2(2u_2+1)+u_1^2)^3
\end{equation*}
and we can spot the equation defining the evolute as the first factor and the outward evolute as the second factor. The real zero locus of this polynomial is plotted in Figure~\ref{discparab}. Observe that the exponent in which appears the outward evolute is $3$, this is related to Corollary~\ref{singdisc}.

\begin{figure}[H]
\centering
\includegraphics[height=5cm, width=5cm]{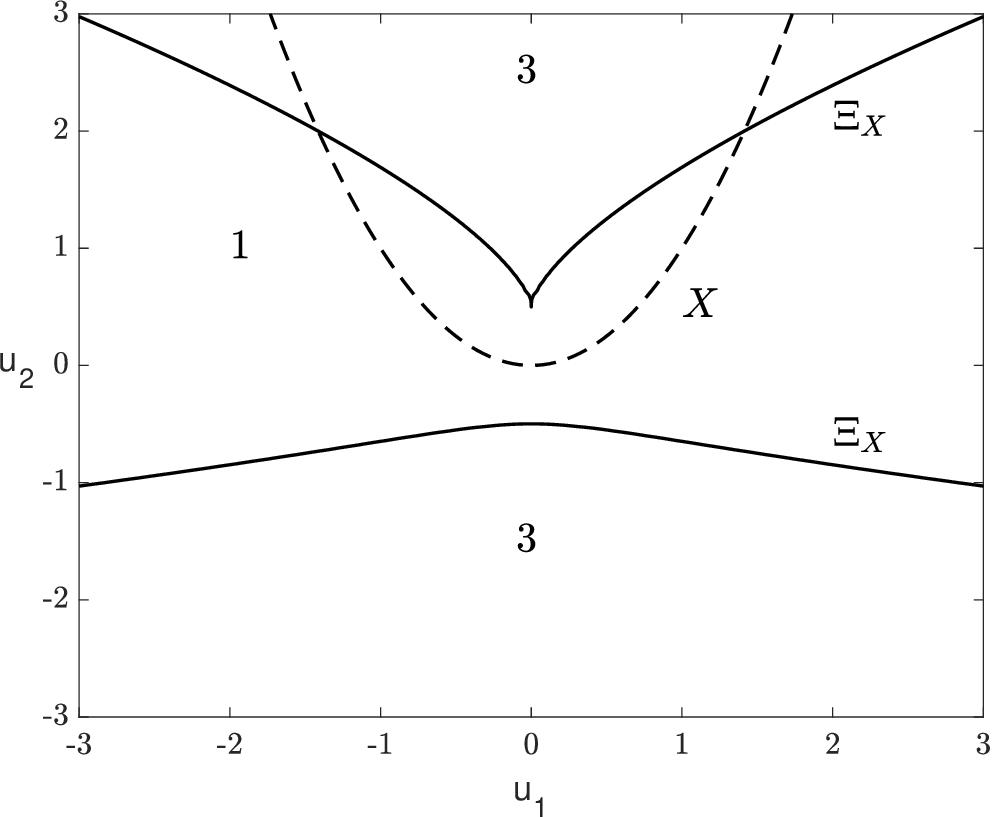}
\caption{Example~\ref{parab1} - Evolute and outward evolute of $X$ (dashed line) on the real plane. The values are the numbers of critical points in each area.}\label{discparab}
\end{figure}

The HD polynomial, see Section~\ref{sec:hdpoly}, of $X$ has degree $5$ in $t^2$ and consist of $131$ monomials. For real points the HD polynomial consists in $79$ terms and decomposes
\begin{align*}
  \HDp_{X,\bu}(t^2)=(4u_2t^2+(u_2-u_1^2)(4u_2+1))^2(16t^6+c_2t^4+c_1t^2+c_0)
\end{align*}
where
\begin{align*}
    &c_2=-8(6u_1^2+2u_2^2+4u_2-1),\\
    &c_1=48u_1^4+32u_1^2u_2^2-8u_1^2u_2+32u_2^3+20u_1^2+8u_2^2-8u_2+1,\\
    &c_0=-(u_1^2-u_2)^2(16u_1^2+(4u_2-1)^2).
\end{align*}
\end{ex}

\begin{ex}\label{circle1}
(Circle). The unit circle 
\begin{equation*}
    X=V(z_1^2+z_2^2-1)\subseteq\C^2
\end{equation*}
satisfies $\vHD(X)=6$ and $\HD=\lbrace 2,4\rbrace$. The value of $\vHD$ can be seen by computing the degree of the critical ideal of $X$ for general parameters $(\bu,\bv)$. This ideal is generated by the polynomials of the system
\begin{equation*}
    \begin{cases}
        p_1(\bz,\bw)=z_1^2+z_2^2-1=0\\
        \bar{p}_1(\bw,\bz)=w_1^2+w_2^2-1=0\\
        p_2(\bz,\bw)=w_2(z_1-u_1)-w_1(z_2-u_2)=0\\ 
        \bar{p}_2(\bw,\bz)=z_2(w_1-v_1)-z_1(w_2-v_2)=0  
    \end{cases}
\end{equation*}
where the upper two polynomials are given by the definition of $X$ and the bottom two polynomials are given by the determinant of the $2\times 2$ matrix obtained from the perpendicularity condition. In particular, it follows the containment $\HD(X)\subseteq\lbrace 2,4,6\rbrace$. We now compute the true values.

Set $\bv=\bar{\bu}$. From the various possible one-to-one parametrizations of the unit circle consider
\begin{align*}
    \psi\colon\C\setminus&\lbrace0\rbrace\to X\subseteq\C^2\\
    &z\mapsto\left(\frac{z^2+1}{2z},\frac{i(z^2-1)}{2z}\right)
\end{align*}
thus critical points of the Hermitian distance satisfy the equation
\begin{align*}
    \scalemath{0.9}{\partial_{z}\|\psi(z)-\bu\|_{2}^2}&\scalemath{0.9}{=\partial_z\left[\left(\frac{z^2+1}{2z}-u_1\right)\left(\frac{\bar{z}^2+1}{2\bar{z}}-\bar{u}_1\right)+\left(\frac{i(z^2-1)}{2z}-u_2\right)\left(\frac{-i(\bar{z}^2-1)}{2\bar{z}}-\bar{u}_2\right)\right]}\\
    &\scalemath{0.9}{=\frac{1}{2z^2}\left((z^2-1)\left(\frac{\bar{z}^2+1}{2\bar{z}}-\bar{u}_1\right)+i(z^2+1)\left(\frac{-i(\bar{z}^2-1)}{2\bar{z}}-\bar{u}_2\right)\right)=0}
\end{align*}
or equivalently, considering the numerator, satisfy the equation
\begin{align*}
    (z^2-&1)\left(\bar{z}^2+1-2\bar{u}_1\bar{z}\right)+(z^2+1)\left(\bar{z}^2-1-2i\bar{u}_2\bar{z}\right)\\
    &=2\left[z^2\bar{z}^2-1-\bar{z}\left((\bar{u}_1+i\bar{u}_2)z^2-\bar{u}_1+i\bar{u}_2\right)\right]=0.
\end{align*}
Since the point $\bu$ is generic, we can solve the complex linear system 
\begin{equation*}
    \begin{cases}
        \bar{a}_1+i\bar{a}_2=-\bar{u}_1\\
        -\bar{a}_1+i\bar{a}_2=-\bar{u}_2
    \end{cases}
\end{equation*}
so that with a change of variables, solving the equation above is equivalent to solving an equation of the form
\begin{equation*}
    p(z,\bar{z})=z^2\bar{z}^2-1+\bar{u}_1z^2\bar{z}+\bar{u}_2\bar{z}=0.
\end{equation*}
By Lemma~\ref{actions}, this problem remains the same under the canonical action on $\C^2$ of the subgroup of the unitary group $\U_2=\U(\C^2)\subseteq\Iso(\C^2)$ that leaves $X$ invariant. This is the group of real orthogonal matrices and thus we can rotate a point $\bu\in\C^2$ to a point $\bu\in\R\times\C$, thus for simplicity we assume $u_1\in\R$. We now consider the system
\begin{equation*}
    \begin{cases}
        zw=\rho\\
        p(z,w)=0\\
        \bar{p}(w,z)=0
    \end{cases}\begin{cases}
        zw=\rho\\
        u_1\rho z+\bar{u}_2w=1-\rho^2\\
        u_2z+u_1\rho w=1-\rho^2
    \end{cases}\quad\cdots\qquad
    \begin{cases}
        zw=\rho\\
        z=\frac{(1-\rho^2)(u_1\rho-\bar{u}_2)}{u_1^2\rho^2-|u_2|^2}\\
        w=\frac{(1-\rho^2)(u_1\rho-u_2)}{u_1^2\rho^2-|u_2|^2}
    \end{cases}
\end{equation*}
for which a solving triple $(z,w,\rho)$ is required to satisfy $\rho>0$ as a necessary condition for the first component $z$ to be a critical point. Thus, we are led to search for positive real solutions of the equation
\begin{equation*}
    zw-\rho=\frac{(1-\rho^2)^2|u_1\rho-\bar{u}_2|^2}{(u_1^2\rho^2-|u_2|^2)^2}-\rho=\frac{u_1^2\rho^6+c_5\rho^5+c_4\rho^4+c_3\rho^3+c_2\rho^2+c_1\rho+|u_2|^2}{(u_1^2\rho^2-|u_2|^2)^2}=0
\end{equation*}
where
\begin{align*}
    &c_5=-u_1(u_1^3+u_2+\bar{u}_2),\\
    &c_4=-2u_1^2+|u_2|^2,\\
    &c_3=2u_1(u_1|u_2|^2+u_2+\bar{u}_2),\\
    &c_2=u_1^2-2|u_2|^2,\\
    &c_1=-|u_2|^4-u_1(u_2+\bar{u}_2).
\end{align*}
Now, we note that for the coefficients $c_4$ to be positive, it requires the coefficient $c_2$ to be negative. To conclude, with this choices the polynomial at the numerator has at most $5$ changes of sign among its ordered coefficients and by the Descartes' rule of signs it admits at most $5$ positive solutions, in particular $\max\HD(X)\leq 5$. By testing points $\bu$ in the real plane we get both two or four critical points and thus it follows $\HD(X)=\lbrace 2,4\rbrace$.

The real locus of the complex evolute of $X$ when $\bv=\bu$ coincides with the union of the evolute and the outward evolute and it is the real zero locus of
\begin{equation*}
    (u_1^2+u_2^2)(u_1^2+u_2^2-4)^3,
\end{equation*}
this locus is plotted in Figure~\ref{disccircle}. Observe that the exponent in which appears the outward evolute is $3$, this is related to Corollary~\ref{singdisc}.

\begin{figure}[H]
\centering
\includegraphics[height=5cm, width=5cm]{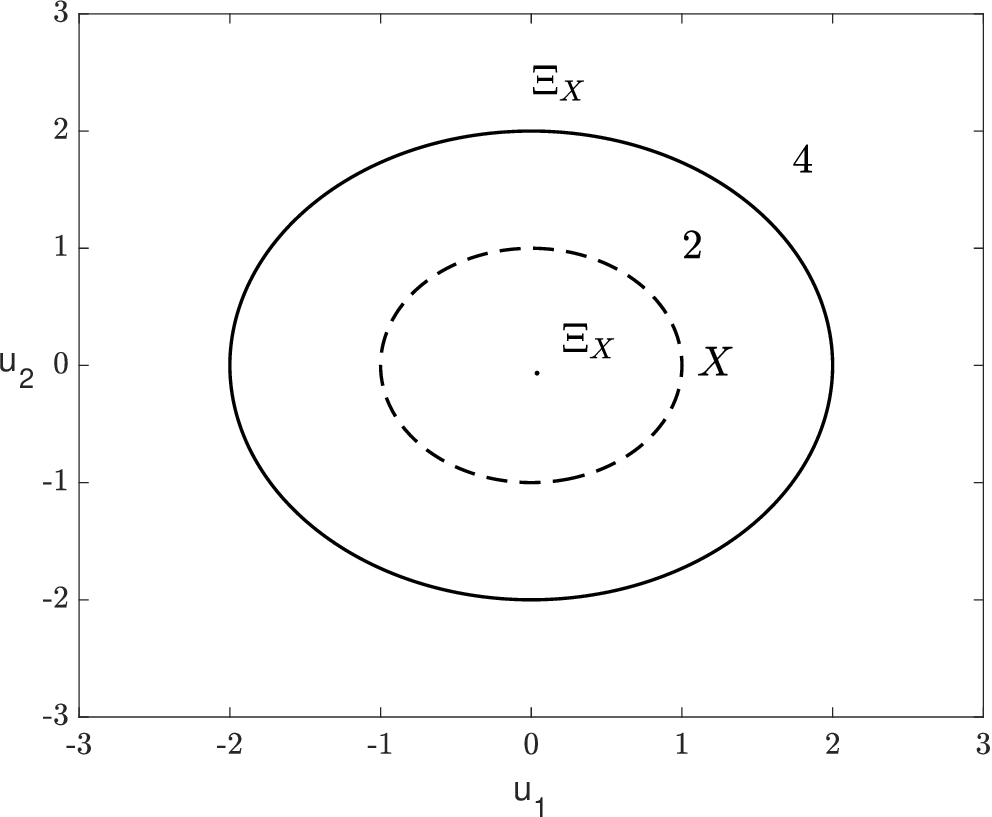}
\caption{Example~\ref{circle1} - Evolute and outward evolute of $X$ (dashed line) on the real plane. The values are the numbers of critical points in each area.}\label{disccircle}
\end{figure}

The HD polynomial, see Section~\ref{sec:hdpoly}, of $X$ has degree $6$ in $t^2$ and consist of $759$ monomials, for real points we get
\begin{equation*}
  \scalemath{0.9}{\HDp_{X,\bu}(t^2)=(t^2-\|\bu\|_{2}^2+1)^2(2t^2-\|\bu\|_{2}^2+2)^2(t^2-2t-\|\bu\|_{2}^2+1)(t^2+2t-\|\bu\|_{2}^2+1)}.
\end{equation*}
\end{ex}

\begin{ex}\label{ellipse1}
(Ellipse). The ellipse 
\begin{equation*}
    X=V(z_1^2+cz_2^2-c)\subseteq\C^2
\end{equation*}
where $0<c<1$ satisfies $\vHD(X)=8$ as predicted by Proposition~\ref{deghyps} and $\lbrace 2,4\rbrace\subseteq\HD(X)\subseteq\lbrace 2,4,6\rbrace$. The value of $\vHD$ can be seen by computing the degree of the Hermitian critical ideal. This ideal is generated by the polynomials of the system
\begin{equation*}
    \begin{cases}
        p_1(\bz,\bw)=z_1^2+cz_2^2-c=0\\
        \bar{p}_1(\bw,\bz)=w_1^2+cw_2^2-c=0\\
        p_2(\bz,\bw)=cw_2(z_1-u_1)-w_1(z_2-u_2)=0\\
        \bar{p}_2(\bw,\bz)=cz_2(w_1-v_1)-z_1(w_2-v_2)=0
    \end{cases}
\end{equation*}
where the upper two polynomials are given by the definition of $X$ and the bottom two polynomials are given by the determinant of the $2\times 2$ matrix obtained from the perpendicularity condition. In particular, it follows the containment $\HD(X)\subseteq\lbrace 2,4,6,8\rbrace$. We now show that $8$ it is not an acceptable value.

Set $\bv=\bar{\bu}$. For the sake of simplicity, we multiply the polynomial defining the variety by $1/c$ and relabel it as $a$ so that $X=V(az_1^2+z_2^2-1)$ where $1<a\in\R$. The matrix representing the Hermitian Killing form
with respect to the basis 
\begin{equation*}
    \lbrace [1],[z_1],[z_2],[w_1],[w_2],[z_2w_2],[w_1w_2],[w_2^2]\rbrace
\end{equation*}
of the quotient $\C[\bz,\bw]/\langle p_1,p_1^{\ast},p_2,p_2^{\ast}\rangle$ is too big to be reported here. However, the determinant of the submatrix given by the subspace $\spn\lbrace [z_1],[z_2]\rbrace$ is 
\begin{equation*}
    -\frac{32a^2(a^2+1)|u_2^2a+u_1^2|^2}{(a^2-1)^4}\leq 0.
\end{equation*}
In particular, this matrix possesses a negative eigenvalue outside the set $V(u_2^2a+u_2^2)\subseteq\C^2$ and from \cite[Corollary 3.14]{f2024} we obtain $\max\HD(X)\leq 7-1=6$.

The real locus of the complex evolute of $X$ when $\bv=\bu$ coincides with the union of the evolute and the outward evolute and it is the zero locus of
\begin{equation*}
  \scalemath{0.9}{\left((4u_1^2+u_2^2-9)^3+972u_1^2u_2^2\right)\left((4u_1^2+u_2^2)^3-880u_1^4-236u_1^2u_2^2+5u_2^4+3900u_1^2-525u_2^2-5625\right)^3},
\end{equation*}
the real zero locus is plotted in Figure~\ref{discEllipse}. Observe that the exponent in which appears the outward evolute is $3$, this is related to Corollary~\ref{singdisc}.

\begin{figure}[H]
\centering
\includegraphics[height=5cm, width=5cm]{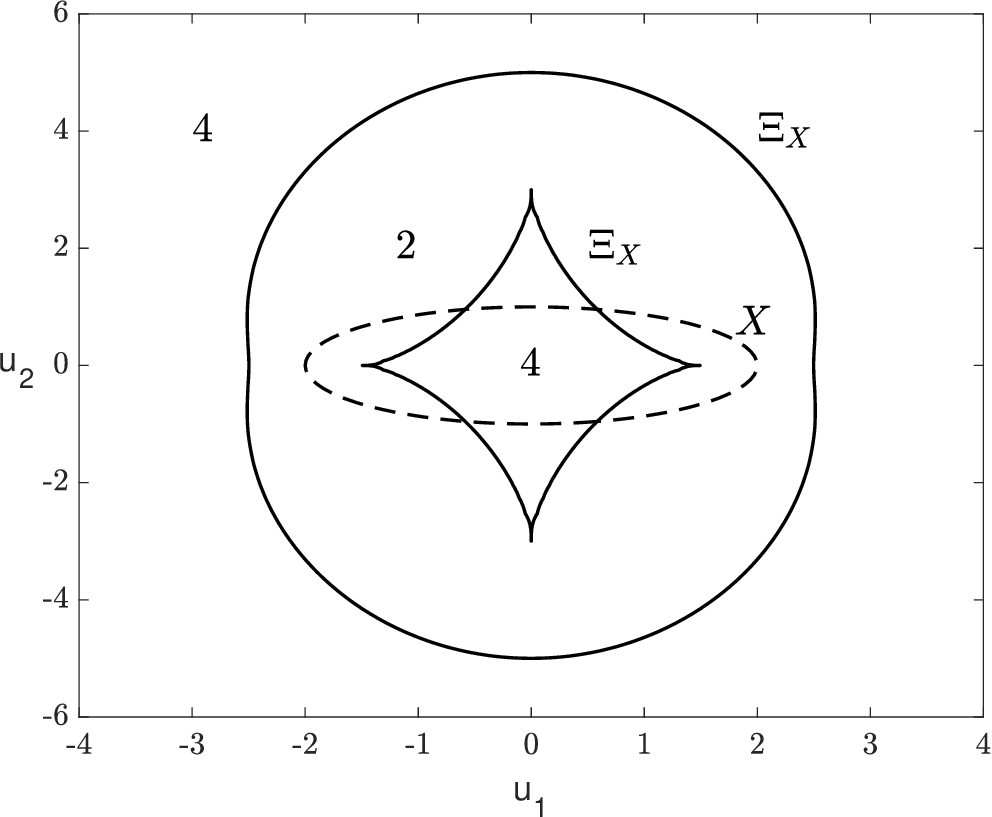}
\caption{Example~\ref{ellipse1} - Evolute and outward evolute of $X$ (dashed line) on the real plane. The values are the numbers of critical points in each area.}\label{discEllipse}
\end{figure}

\end{ex}


\section{The HD polynomial}\label{sec:hdpoly}\thispagestyle{plain}

In this section we introduce a polynomial that encodes all the notions we introduced so far about the Hermitian distance problem, from the exact distances from the critical points to the HD discriminant of Section~\ref{sec:hddisc}. The \emph{Euclidean distance polynomial} (ED polynomial) introduced in \cite{os2020} is the equivalent concept in the case of the Euclidean distance problem, a consistent study about it can be found in \cite{hw2018}, the same concept is called the \emph{offset polynomial} in \cite{bks2024}.

\begin{Def}
Let $X$ be an algebraic variety. The sum of the Hermitian critical ideal $\eqref{cond}$ and the ideal $\langle t^2-\langle\bu-\bz,\bv-\bw\rangle_{\R}\rangle$ in the ring $\C[\bz,\bw,\bu,\bv,t]$ defines a variety of dimension $2n$.

The \emph{virtual Hermitian Distance polynomial} (vHD polynomial) $\vHDp_{X,\bu,\bv}(t^2)$ of $X$ at $(\bu,\bv)$ is, up to a scalar factor, the generator of the projection of the sum of ideals above in the ring $\C[\bu,\bv,t]$. When the variety $X$ is defined by polynomials with real (rational) coefficients, the elimination procedure yields $\vHDp_{X,\bu,\bv}$ with real (rational) coefficients. Moreover, the \emph{Hermitian Distance polynomial} (HD polynomial) is $\HDp_{X,\bu}\coloneqq\vHDp_{X,\bu,\bar{\bu}}$.
\end{Def}

Since the polynomial ring is a UFD, the projection of the sum of the ideals is generated by a single polynomial in $t^2$ and thus this object is well defined.\newline

Recall from Section~\ref{sec:first} that we denote $X_{\circ}=X_{\reg}\times\overline{X}_{\reg}$. From the definition of vHD polynomial we obtain the following result.

\begin{Prop}\label{hdp}
Let $X$ be an algebraic variety and $(\bu,\bv)$ be a generic point. The roots of the polynomial $\vHDp_{X,\bu,\bv}(t^2)$ are of the form $t^2=\langle\bu-\bz,\bv-\bw\rangle_{\R}$, where $(\bz,\bw)\in X_{\circ}$ is a solution of the Hermitian critical ideal of $(\bu,\bv)$. In particular, the distance of $\bu$ from $X$ is a root of $\HDp_{X,\bu}(t^2)$ and the degree of $\vHDp_{X,\bu,\bv}(t^2)$ in $t^2$ is equal to $\vHD(X)$.
\end{Prop}

The next proposition highlights a useful property of the HD polynomial.

\begin{Prop}
Consider the pair-wise disjoint union $X=X_1\cup\ldots\cup X_r$ for some $r\in\N$, where $X_k\subseteq\V$ is a reduced variety for every $k=1,\ldots,r$. It holds the equality 
\begin{equation*}
    \vHDp_{X,\bu,\bv}(t^2)=\prod_{k=1}^r\vHDp_{X_k,\bu,\bv}(t^2).
\end{equation*}
\end{Prop}
\begin{proof}
The variety of the Hermitian critical ideal \eqref{cond} of $X$ is the union of the varieties of the Hermitian critical ideals \eqref{cond} of the varieties $X_k$. The conclusion follows by Proposition~\ref{hdp}.
\end{proof}

The following lemma generalizes Remark~\ref{diltras}.

\begin{Lemma}
Let $X\subseteq\V$ be an algebraic variety, $0\neq c\in\C$ and $\bb\in\V$, it holds the equality
\begin{equation*}
\vHDp_{cX+\bb,c\bu+\bb,\bar{c}\bv+\bar{\bb}}(|c|^2t^2)=\vHDp_{X,\bu,\bv}(t^2).
\end{equation*}
\end{Lemma} 
\begin{proof}
The assertion follows from Remark~\ref{diltras} by noting 
\begin{equation*}
    \langle c\bu+\bb-(c\bz+\bb),\bar{c}\bv+\bar{\bb}-(\bar{c}\bw+\bar{\bb})\rangle_{\R}=|c|^2\langle\bu-\bz,\bv-\bw\rangle_{\R}.
\end{equation*}
\end{proof}

We denote the discriminant of the vHD polynomial as $\Delta_{X}\coloneqq\Delta_{t^2}\vHDp_{X,\bu,\bv}(t^2)$. We get the following result connecting the vHD polynomial and the vHD discriminant.

\begin{Prop}
Let $X$ be an algebraic variety. The polynomial defining the vHD discriminant divides $\Delta_{X}$. Moreover, if $X$ is symmetric with respect to $r$ affine hyperplanes $L_1,\ldots,L_r$ of equations $l_1,\ldots,l_r$, then the product $l_1\cdots l_rl_1^{\ast}\cdots l_r^{\ast}$ divides $\Delta_{X}$, where $\ast$ is the map of Section~\ref{sec:first} and we are assuming $l_1,\ldots,l_r$ to be polynomials in $\bz$.
\end{Prop}
\begin{proof}
By definition of $\vX_{X}$, any point $(\bu,\bv)\in\vX_{X}$ satisfies the equation of $\Delta_{X}$, since two roots in $t^2$ coincide by Proposition~\ref{hdp}. Now, let $(\bu,\bv)\in L_k\times\overline{L}_j$ be generic and let 
\begin{equation*}
    (\bz,\bw)\in \left(X\setminus(X_{\sing}\cup L_k)\right)\times\left(\overline{X}\setminus(\overline{X}_{\sing}\cup\overline{L}_{j})\right)
\end{equation*}
for some $k,j\in\lbrace 1,\ldots,r\rbrace$ be a solution of the Hermitian critical ideal of $(\bu,\bv)$. Denote by $(\tilde{\bz},\tilde{\bw})$ the reflection of $(\bz,\bw)$ with respect to $L_k\times\overline{L}_j$. Thus $(\tilde{\bz},\tilde{\bw})\in X_{\circ}$ as well and $(\tilde{\bz},\tilde{\bw})$ is again a solution of the Hermitian critical ideal of $(\bu,\bv)$. Since $\langle\bu-\bz,\bv-\bw\rangle_{\R}=\langle\bu-\tilde{\bz},\bv-\tilde{\bw}\rangle_{\R}$ then $(\bu,\bv)$ is a zero of $\Delta_{X}$.
\end{proof}

\begin{Prop}
Let $X\subseteq\V$ be a real algebraic variety for which it holds $\max\HD(X)<\vHD(X)$. The set $\lbrace(\bu,\bu)\in\V^2\mid\bu\in\V\rbrace\subseteq\V^2$ is in the zero locus of $\Delta_{X}$.
\end{Prop}
\begin{proof}
Similarly as discussed in Remark~\ref{realpoints}, by the hypothesis there exists at least a solution $(\bz,\bw)$ of the Hermitian critical ideal of $(\bu,\bu)\in\V^2$ such that $\bz\neq\bw$ which in turn yields another solution $(\bw,\bz)$ and those solutions satisfy
$\langle\bu-\bz,\bu-\bw\rangle_{\R}=\langle\bu-\bw,\bu-\bz\rangle_{\R}=t^2$. In particular, the HD polynomial has a solution with multiplicity greater than $1$ and the assertion follows by definition of $\Delta_{X}$.
\end{proof}

The proposition above highlights the fact that for a real point there exist two critical points at the same distance which are conjugate one to each other.

\begin{Prop}
Let $X\subseteq\V$ be an algebraic variety and $G\subseteq\Iso(\V)$ be a subgroup that leaves $X$ invariant. For any $g\in G$ it holds $\vHDp_{X,\bu,\bv}=\vHDp_{X,g\cdot\bu,\bar{g}\cdot\bv}$ and the coefficients of $\HDp_{X,\bu}$ are $G$ invariant.
\end{Prop}
\begin{proof}
Follows from Lemma~\ref{actions}.
\end{proof}

The following result is the Hermitian version of \cite[Proposition 2.12]{os2020}.

\begin{Prop}\label{affsub}
Let $X$ be an affine subspace, there hold $\vHD(X)=1$ and $\HD(X)=\lbrace 1\rbrace$. 
Moreover, if $\bb\in X$ then 
\begin{equation*}
	\vHDp_{X,\bu,\bv}(t^2)=t^2-\langle\pi_{\overline{(X-\bb)}^{\perp_{\R}}}(\bu-\bb),\pi_{(X-\bb)^{\perp_{\R}}}(\bv-\bar{\bb})\rangle_{\R}.
\end{equation*}
\end{Prop}
\begin{proof}
The tangent spaces of $X$ and $\overline{X}$ are everywhere isomorphic to $\spn\lbrace e_k\rbrace_{k=1}^m$ and $\spn\lbrace \bar{e}_k\rbrace_{k=1}^m$ respectively, where $m\in\N$ is the dimension of $X$ and $e_1,\ldots,e_m\in\V$ are orthonormal generators with respect to the Hermitian inner product. If the point $(\bz,\bw)$ is a solution of the Hermitian critical ideal of $(\bu,\bv)\in\V^2$, then $\bv-\bw\perp_{\R}\spn\lbrace e_k\rbrace_{k=1}^m$ and $\bu-\bz\perp_{\R}\spn\lbrace \bar{e}_k\rbrace_{k=1}^m$. We complete $\lbrace e_k\rbrace_{k=1}^m$ with vectors $e_{m+1},\ldots,e_{n}$ to get an orthonormal basis of $\V$ with respect to the Hermitian inner product. By writing $\bv=v_1\bar{e}_1+\ldots+v_n\bar{e}_n$ and $\bu=u_1e_1+\ldots+u_ne_n$, it is easy to check that the only possible solution is 
\begin{equation*}
	\bz=(u_1-\langle\bb,\bar{e}_1\rangle_{\R})e_1+\ldots+(u_m-\langle\bb,\bar{e}_m\rangle_{\R})e_m+\bb\in X\quad\text{and}\quad\bw=\bar{\bz}\in\overline{X}.
\end{equation*}
The equation for the HD polynomial follows.

Clearly, the only critical point is $\pi_{(X-\bb)}(\bu-\bb)+\bb$. 
\end{proof}

\begin{Coro}
Let $X=V(a_1z_1+\ldots+a_nz_n+b)$ with $a_1,\ldots,a_n\in\C$ not all zero and $b\in\C$ be a hyperplane and denote $\ba=(a_1,\ldots,a_n)$. It holds 
\begin{equation*}
	\vHDp_{X,\bu,\bv}(t^2)=t^2-\frac{\langle\|\ba\|_{2}^2\bu-b\bar{\ba},\ba\rangle_{\R}\langle\|\ba\|_{2}^2\bv-\bar{b}\ba,\bar{\ba}\rangle_{\R}}{\|\ba\|_{2}^4}.
\end{equation*} 
\end{Coro}

\begin{ex}\label{exlin}
(Linear regression). We can bijectively parametrize a hyperplane as
\begin{align*}
&\psi\colon\C^{n-1}\longrightarrow X=V\left(z_n-\sum_{k=1}^{n-1}a_kz_k-b\right)\subseteq\C^n\\
&(z_1,\ldots,z_{n-1})\mapsto\left(z_1,\ldots,z_{n-1},\sum_{k=1}^{n-1}a_kz_k+b\right)
\end{align*}
Critical points of $q_{\bu}$ satisfy the following equation 
\begin{equation*}
   \partial_{z_k}\|\psi(z_1,\ldots,z_{n-1})-\bu\|_{2}^2=(\bar{z}_k-\bar{u}_k)+a_k \left(\sum_{k=1}^{n-1}\bar{a}_k\bar{z}_k+\bar{b}-\bar{u}_n\right)=0
\end{equation*}
for any $k=1,\ldots,n-1$. Denoting the row vectors $\ba=(a_1,\ldots,a_{n-1})$, $\tilde{\bu}=(u_1,\ldots,u_{n-1})$ and $\tilde{\bz}=(z_1,\ldots,z_{n-1})$, in other terms we aim to solve the $(n-1)\times(n-1)$ system
\begin{equation*}
 \left(I_{n-1}+\ba^H\ba\right)\tilde{\bz}^T=(u_n-b)\ba^H+\tilde{\bu}^T.
\end{equation*}
We can apply the inverse $I_{n-1}-\frac{1}{1+\|\ba\|_{2}^2}\ba^H\ba\in\C^{(n-1)\times(n-1)}$ of the matrix on the left hand side to find $\tilde{\bz}$. Thus, by computing the last coordinate $z_n=\langle \ba,\tilde{\bz}\rangle_{\R}+b$ we obtain the explicit unique critical point 
\begin{equation*}
(\tilde{\bz},z_n)=\left(\frac{u_n-b-\langle\ba,\tilde{\bu}\rangle_{\R}}{1+\|\ba\|_{2}^2}\bar{\ba}+\tilde{\bu},\frac{u_n\|\ba\|_{2}^2+b+\langle\ba,\tilde{\bu}\rangle_{\R}}{1+\|\ba\|_{2}^2}\right)\in X.
\end{equation*}
We parametrize the HD correspondence as
\begin{align*}
  &\C^{n-1}\times\C\longrightarrow\mathcal{H}_{X}\subseteq\C^n\times\C^n\\
  (\tilde{\bz}&,s)\mapsto\left((\tilde{\bz},\langle\ba,\tilde{\bz}\rangle_{\R}+b),(\tilde{\bz}+s\bar{\ba},\langle\ba,\tilde{\bz}\rangle_{\R}+b-s)\right).
\end{align*}
If $X=V(z_1-a_1,\ldots,z_{n-1}-a_{n-1})\subseteq\C^n$, then the critical point is $(\ba,u_n)\in X$.
\end{ex}

The following result relates the HD polynomial of a variety $X$ to the one of its dual variety $X^{\vee}$, see Subsection~\ref{ssec:corrdual}.

\begin{Theo}\label{poldual}
Let $X$ be an affine cone, it holds
\begin{equation*}
    \vHDp_{X,\bu,\bv}(t^2)=\HDp_{X^{\vee},\bv,\bu}(\langle\bu,\bv\rangle_{\R}-t^2).
\end{equation*}
\end{Theo}
\begin{proof}
Let $(\bz,\bw)\in X_{\circ}$ be a zero of the Hermitian critical ideal of $X$ of $(\bu,\bv)\in\V^2$. By Proposition~\ref{hdp} the value $\langle\bu-\bz,\bv-\bw\rangle_{\R}$ is a root of $\vHDp_{X,\bu,\bv}$. Applying Theorem~\ref{dual} the point $(\bv-\bw,\bu-\bz)$ is a solution of the Hermitian critical ideal of $X^{\vee}$ of $(\bv,\bu)$. Hence, again by Proposition~\ref{hdp}, the value $\langle\bw,\bz\rangle_{\R}$ is a root of $\HDp_{X^{\vee},(\bv,\bu)}$. In the end, since $X$ is an affine cone, summing the two roots we get
\begin{equation*}
    \langle\bu-\bz,\bv-\bw\rangle_{\R}+\langle\bw,\bz\rangle_{\R}=\langle\bu,\bv\rangle_{\R}-\langle\bz,\bv-\bw\rangle_{\R}-\langle\bu-\bz,\bw\rangle_{\R}=\langle\bu,\bv\rangle_{\R}
\end{equation*}
and thus we can rewrite the zeros of $\HDp_{X^{\vee},(\bv,\bu)}$ as
\begin{equation*}
   \langle\bw,\bz\rangle_{\R}=\langle\bu,\bv\rangle_{\R}-\langle\bu-\bz,\bv-\bw\rangle_{\R}.
\end{equation*}
\end{proof}

The next result states the obvious fact that the points of the variety are at distance zero from the variety itself.

\begin{Prop}
The algebraic variety $X$ is contained in the zero locus of $\HDp_{X,\bu}(0)$.
\end{Prop}

\begin{ex}\label{exfcurve}
(Fermat curve of degree $2$). The Fermat curve 
\begin{equation*}
    X=V(z_1^2+z_2^2-z_3^2)\subseteq\Pro^2
\end{equation*}
is the projective closure of the variety of Example~\ref{circle1}, moreover it satisfies the equality $X^{\vee}=X$. From Theorem~\ref{projcon}, we already know $\vHD(X)=2$ and $\HD(X)=\lbrace 2\rbrace$. The polynomial defining the vHD discriminant $\vX_{X}$ is
\begin{equation*}
    ((v_1u_2-u_1v_2)^2-(v_1u_3+u_1v_3)^2-(v_2u_3+u_2v_3)^2)^2
\end{equation*}
and in particular the generalized polynomial defining the HD discriminant $\Xi_{X}$ is
\begin{equation*}
    ((\bar{u}_1u_2-u_1\bar{u}_2)^2-(\bar{u}_1u_3+u_1\bar{u}_3)^2-(\bar{u}_2u_3+u_2\bar{u}_3)^2)^2.
\end{equation*}
The vHD polynomial is
\begin{equation*}
    \vHDp_{X,\bu,\bv}(t^2)=4t^4-4(u_1v_1+u_2v_2+u_3v_3)t^2+(u_1^2+u_2^2-u_3^2)(v_1^2+v_2^2-v_3^2),
\end{equation*}
in particular we obtain
\begin{equation*}
    \HDp_{X,\bu}(t^2)=4t^4-4\|\bu\|_{2}^2t^2+|u_1^2+u_2^2-u_3^2|^2.
\end{equation*}
\end{ex}


\section{Conclusion}\label{sec13}

We conclude this investigation of the Hermitian Distance degree by adapting the conclusion from \cite{dhost2014}:

\begin{centering}
    \textit{The simplest curves are the planar curves. Among them, the simplest one is the line $(\HD=\lbrace 1\rbrace)$. The next simplest curve is the parabola $(\HD=\lbrace 1,3\rbrace)$. After that come the circle $(\HD=\lbrace 2,4\rbrace)$.}
\end{centering}

Unfortunately, the case of the ellipse remains to be solved.

\backmatter


\bmhead{Acknowledgements}

I would like to thank Prof.\ Giorgio Ottaviani who propose the problem to me and supervised me throughout the redaction of this work.

\bibliography{sn-bibliography}%

@article{dhost2014,
  title     = {The {E}uclidean {D}istance {D}egree of an {A}lgebraic {V}ariety},
  author    = {Jan Draisma and Emil Horobe\c{t} and Giorgio Ottaviani and Bernd Sturmfels and Rekha R Thomas},
  journal   = {Foundations of Computational Mathematics},
  volume    = {16},
  pages     = {99--149},
  year      = {2016}
}

@article{f2024,
  title     = {The {H}ermitian {K}illing form and root counting of complex polynomials with conjugate variables},
  author    = {Davide Furch\`i},
  journal   = {Linear Algebra and its Applications},
  volume    = {708},
  pages     = {93--111},
  year      = {2025}
}

@article{o2022,
  title     = {The critical space for orthogonally invariant varieties},
  author    = {Giorgio Ottaviani},
  journal   = {Vietnam Journal of Mathematics},
  volume    = {50},
  pages     = {615-622},
  year      = {2022}
}

@article{c2019,
  title     = {Unmixing the mixed volume computation},
  author    = {Tianran Chen},
  journal   = {Discrete \& Computational Geometry},
  volume    = {62},
  pages     = {55--86},
  year      = {2019}
}

@book{e1996,
  author    = {G\"unter Ewald},
  year      = {1996},
  title     = {Combinatorial {C}onvexity and {A}lgebraic {G}eometry},
  publisher = {Springer},
  address   = {New York}
}

@book{hj2013,
  author    = {Roger A Horn and Charles R Johnson},
  year      = {2013},
  title     = {Matrix {A}nalysis},
  publisher = {Cambridge University Press},
  address   = {New York}
}

@book{gkz1994,
  author    = {Israel M Gelfand and Mikahil M Kapranov and Andrei V Zelevinsky},
  year      = {1994},
  title     = {Discriminants, resultants and higher-dimensional determinants},
  publisher = {Birkh\"auser},
  address   = {Boston}
}

@book{bg1992,
  author    = {James W Bruce and Peter Giblin},
  year      = {1992},
  title     = {Curves and {S}ingularities},
  publisher = {Cambridge University Press},
  address   = {New York}
}

@article{c1998,
  title     = {Focal {L}oci of {A}lgebraic {H}ypersurfaces: {A} {G}eneral {T}heory},
  author    = {Cecilia Trifogli},
  journal   = {Geometriae Dedicata},
  volume    = {70},
  pages     = {1--26},
  year      = {1998}
}

@article{ct2007,
  title     = {Focal loci of algebraic varieties {I}},
  author    = {Fabrizio Catanese and Cecilia Trifogli},
  journal   = {Communications in Algebra},
  volume    = {28},
  number    = {12},
  pages     = {6017--6057},
  year      = {2007}
}

@article{os2020,
  title     = {The distance function from a real algebraic variety},
  author    = {Giorgio Ottaviani and Luca Sodomaco},
  journal   = {Computer Aided Geometric Design},
  volume    = {82},
  pages     = {101927},
  year      = {2020}
}

@article{hw2018,
  title     = {Offset {H}ypersurfaces and {P}ersistent {H}omology of {A}lgebraic {V}arieties},
  author    = {Emil Horobe\c{t} and Madeleine Weinstein},
  journal   = {Computer Aided Geometric Design},
  volume    = {74},
  number    = {3},
  pages     = {1521--1542},
  year      = {2018}
}

@book{bks2024,
	author    = {Paul Breiding and Kathl\'{e}n Kohn and Bernd Sturmfels},
	year      = {2024},
	title     = {Metric {A}lgebraic {G}eometry},
	publisher = {Birkh\"{a}user},
    address   = {Cham}
}

\end{document}